\magnification 1200
\input amssym.def
\input amssym.tex
\parindent = 40 pt
\parskip = 12 pt
\font \heading = cmbx10 at 12 true pt
 at 22 true pt
 at 19 true pt
 at 7 true pt
\def \R{{\bf R}}

\centerline{\heading Resolution of Singularities in Two Dimensions}
\centerline{\heading and the Stability of Integrals}
\rm
\line{}
\line{}
\centerline{\heading Michael Greenblatt}
\line{}
\centerline{June 9, 2009}
\baselineskip = 12 pt
\font \heading = cmbx10 at 14 true pt
\line{}
\line{}
\noindent{\bf 1. Background and Statement of Results}

\vfootnote{}{This research was supported in part by NSF grant DMS-0654073} Suppose $S(x,y)$ is a smooth
real-valued function on a neighborhood of $(0,0)$ with $S(0,0) = 0$. For a small open set $U$ containing 
the origin and a small $\epsilon > 0$, define $M_{S,U}(\epsilon)$ by
$$M_{S,U}(\epsilon) = |\{(x,y) \in U: |S(x,y)| < \epsilon\}| \eqno (1.1)$$
By resolution of singularities (or using the explicit formulas of [G1]) if $S(x,y)$ is real-analytic then 
if $U$ is sufficiently small then for some $C_U > 0$ one has
$$M_{S,U}(\epsilon) = C_U \epsilon^j |\ln \epsilon|^p + o(\epsilon^j |\ln\epsilon|^p)\eqno (1.2)$$
Here $j$ is a positive rational number and $p = 0$ or 1. Also using resolution of singularities, it can be
shown that $(j,p)$ is independent of $U$ for small enough $U$. We refer to $j$ as the 
{\it growth index} of $S$ and to $p$ as the {\it multiplicity} of this growth index. If $S(x,y)$ is merely smooth,
then it follows from [G1] that $(1.2)$ still holds unless there is a smooth coordinate change fixing the 
origin after which the bisectrix intersects the Newton polygon of $S$ in the interior of its horizontal
ray (see below for the relevant definitions.) It is also true that except in those exceptional 
situations, $C_U$ is independent of $U$ for small enough $U$. In the exceptional situations, then  
we can at least say there is a $j$ such that for small enough $U$, for some $C_U > 0$ we have 
$M_{S,U}(\epsilon) < C_U\epsilon^j$, while for $j' > j$ there is no $C_U'$ for which the estimate 
$M_{S,U}(\epsilon) < C_U'\epsilon^{j'}$ holds. Thus for the smooth situation, we also have a natural 
definition of the growth index of $S(x,y)$ and its multiplicity. 

A natural question to consider
is the effect of perturbing $S(x,y)$ on the growth index and its multiplicity. Besides being of 
intrinsic interest, these questions and their oscillatory integral analogues (described later in this
section) are important 
in the analysis of Fourier transforms of surface-supported measures such as in [IKeM] and [IoSa].
Complex and higher-dimensional analogues of these questions are also connected to various issues in 
complex geometry. 

Since by the implicit 
function theorem one has $M_{S,U}(\epsilon) \sim |\nabla S (0,0)|^{-1}\epsilon$ when $\nabla S(0,0) \neq 
0$ and $U$ is sufficiently small, one typically assumes that $S(x,y)$ has a critical point at the origin.
That is, one assumes that $S(0,0) = 0$ and $\nabla S(0,0) = 0$.

For the analogous problem in one dimension, the growth index is just the reciprocal of the order of 
vanishing of $S(x)$ at $x = 0$, so a small perturbation of the phase can only result in the growth index 
staying the same or increasing. A famous example of Varchenko [V] shows that the analogous phenomenon 
does not necessarily hold in three or higher dimensions. Thus a general result can hold only in dimension two.
The strongest result in the current literature is the following theorem due to Karpushkin [K3]. Let $D_r$ 
denote the open disk in $\R^2$ of radius $r$ centered at the origin, and $E_r$ the open disk in 
${\bf C}^2$ of radius $r$ centered at the origin. For a function $f(x,y)$ real-analytic on $D_r$, let
$\tilde{f}(z_1,z_2)$ denote the unique holomorphic extension of $f(x,y)$ to $E_r$. Then Karpushkin's 
theorem is:

\noindent {\bf Theorem:} ([K3]) Suppose $S(x,y)$ is real-analytic on $D_r$ satisfying $S(0,0) = 0$ and
$\nabla S(0,0) = 0$ such that $S(x,y)$ has 
growth index $j$ with multiplicity $p$ at the origin. Then there is a $\delta > 0$, an $s < r$, 
and a positive constant $C_S$ depending on $S(x,y)$ such that if $f(x,y)$ is real-analytic on $D_s$
and $\tilde{f}(z_1,z_2)$ extends to a continuous function on $\bar{E}_s$ with $|\tilde{f}(z_1,z_2)| 
< \delta$ for all $(z_1,z_2) \in \bar{E}_s$ then for $0 < \epsilon < {1 \over 2}$ 
$$M_{S + f, D_s}(\epsilon) \leq C_S \epsilon^j |\ln\epsilon|^p \eqno (1.3)$$
Thus Karpushkin's theorem shows that for real-analytic functions not only does the growth index and 
multiplicity either improve or stay the same under a small perturbation, but also one has uniformity in 
the radius $s$ and in the constant $C_S$. Karpushkin's proofs involve ideas from singularity
theory, in particular the theory of versal deformations which turns arbitrary perturbations of $S(x,y)$
into a number of canonical forms which then may be considered individually. 

Another method for dealing with stability of $M_{S,U}(\epsilon)$ was introduced in [PSSt] where the
above theorem is proven modulo the logarithmic factors. Their methods are often referred to in this 
subject as the method 
of algebraic estimates, which also give partial analogues in higher dimensions (the example of Varchenko
shows the full analogues not feasible). Also, in the case of {\it linear} perturbations of smooth functions,
results of the above nature are proven in [IKeM]. 

The purpose of this paper is to show how resolution of singularities algorithms in two
dimensions, in conjunction with some one-dimensional Van der Corput-type lemmas, provides another method 
which we will use to prove new estimates and 
theorems for the $M_{S,U}(\epsilon)$ as well as for oscillatory integral analogues. Since these
algorithms will apply to all smooth functions, our theorems will hold for all smooth functions as opposed to the 
earlier real-analytic results of [K1]-[K3] and [PSSt]. We also will make use of the superadapted coordinate 
systems of [G1] that put functions in certain canonical forms suitable for these problems. They are a
refinement of the adapted coordinate systems of [V]. Adapted coordinate systems are also used in [K1]-[K3]
and [PSSt].

For our sharpest estimates (Theorems 1.1 and 1.5) our condition on the perturbation function $f(x,y)$ 
will be that the absolute values of finitely many derivatives of $f(x,y)$ at the origin are less than some $\delta$ which
depends on $S(x,y)$. We will get uniformity in the coefficient $C_S$ of $(1.3)$, but there will not be 
uniformity in the radius $s$ since such uniformity does not hold in the general smooth case. For a given
$f(x,y)$, there will also be finitely many $t$ with $|t| < 1$ such that the case when 
$tf(x,y)$ is the perturbation function is excluded from the theorems. This is due to certain error terms 
being affected by the zeroes of certain one-dimensional polynomials induced by $S(x,y) + tf(x,y)$. Since 
these issues only affect error terms, other than for these exceptional values we will get the uniform 
estimates. 

\noindent Before we state our theorems we first define some relevant terminology.

\noindent {\bf Definition 1.1.} Let $S(x,y) = \sum_{a,b} s_{ab}x^ay^b$ denote the 
Taylor expansion of $S(x,y)$ at the origin. Assume there is at least one $(a,b)$ for which $s_{ab}$ is
nonzero. For any $(a,b)$ for which $s_{ab} \neq 0$, let $Q_{ab}$ be the quadrant $\{(x,y) \in \R^2: 
x \geq a, y \geq b \}$. Then the {\it Newton polygon} $N(S)$ of $S(x,y)$ is defined to be 
the convex hull of the union of all $Q_{ab}$.  

A Newton polygon consists of finitely many (possibly zero) bounded edges of negative slope
as well as an unbounded vertical ray and an unbounded horizontal ray. More generally, one can define the 
Newton polygon of a power series in $x^{1 \over N}$ and $y^{1 \over N}$ for a positive integer $N$ 
analogously to Definition 1.1.

\noindent {\bf Definition 1.2.} The {\it Newton distance} $d(S)$ of $S(x,y)$ is defined to be 
$\inf \{t: (t,t) \in N(S)\}$.

Throughout this paper, we will use the $(a,b)$ coordinates to write equations of lines relating
to Newton polygons, so as to distinguish from the $x$-$y$ variables of the domain of $S(x,y)$. The line 
in the $a$-$b$ plane with equation $a = b$ comes up so frequently it has its own name:

\noindent {\bf Definition 1.3.} The {\it bisectrix} is the line in the $a$-$b$ plane with equation
$a = b$.

\noindent In Theorems 1.1-1.3 below, $S(x,y)$ is a smooth function on a neighborhood of the origin 
with nonvanishing Taylor expansion at the origin and satisfying $S(0,0) = 0$ and $\nabla S(0,0) = 0$. 
Denote the growth index of $S(x,y)$ by $j$ and its multiplicity by $p$. Our first and sharpest theorem
is the following.

\noindent {\bf Theorem 1.1.} There is a positive integer $l$ and a $\delta > 0$ such that if $f(x,y)$ 
is a smooth function on a neighborhood of the origin with $\sup_{|\alpha| \leq l} |\partial^{\alpha}f
(0,0)| < \delta$, then for all but finitely many $t$ with $|t| < 1$, if $D$ is a sufficiently small disk 
centered at the origin (depending on $S + tf$), for all $0 < \epsilon < {1 \over 2}$ we have 
$$M_{S + tf,D}(\epsilon) \leq C_S \epsilon^j |\ln \epsilon|^p \eqno (1.4)$$

There is no uniformity in the radius of $D$ in Theorem 1.1 as such a statement is false for
general smooth functions. It should also be pointed out that since in most cases the leading coefficient of 
$(1.2)$ is independent of $U$ for small enough $U$, simply shrinking down $D$ does not typically help
in getting a uniform constant in the right-hand side of $(1.4)$.

Our next result says that as long as $S(x,y)$ and 
$f(x,y)$ do not both have Morse (nondegenerate) critical points at the origin, then Theorem 1.1 still holds for $|t| \geq 1$ 
as well, although now the constant now depends on $t$ and $f$ as well as $S$. Another way of saying this is
that the growth index and multiplicity of $S(x,y) + tf(x,y)$ is at least as good as that of $S(x,y)$ for all
but finitely many $t$.  

\noindent {\bf Theorem 1.2.} Let $f(x,y)$ be a smooth function on a neighborhood of the origin with
a critical point there and assume that $S(x,y)$ and $f(x,y)$ do not both have Morse critical 
points at the origin. Let $j_t$ denote the growth index of $S(x,y) + tf(x,y)$ at the origin and $p_t$ 
its multiplicity. Then for all but finitely
many real values of $t$ we have $(-j_t,p_t) \leq (-j,p)$ (under the lexicographic ordering). 

Theorem 1.2 does not hold if $S(x,y)$ and $f(x,y)$ are both Morse as can be seen by taking
$S(x,y) = x^2 + y^2$ and $f(x,y) = x^2 - y^2$. However, in such situations the growth index of all but
finitely many $S + tf$ is still going to be $1$. Also, note that the condition excluding finitely many $t$
may be necessary; for example, when $f(x,y) = -S(x,y)$ plus a small error term. 

Since the growth index and multiplicity of $\alpha S_1(x,y) + \beta S_2(x,y)$ is the same as 
that of $S_1(x,y) + {\beta \over \alpha} S_2(x,y)$ for any $\alpha \neq 0$, Theorem 1.2 and symmetry 
imply the following:

\noindent {\bf Theorem 1.3.} Suppose $S_1(x,y)$ and $S_2(x,y)$ are two smooth functions on a neighborhood
of the origin with critical points at the origin, not both Morse. Let $(j_1,p_1)$ and $(j_2,p_2)$ be their growth 
indices and multiplicities. For any (real) $\alpha$ and $\beta$, Let $j_{\alpha,\beta}$ and $p_{\alpha,\beta}$ be the 
growth index and multiplicity of $\alpha S_1(x,y) + \beta S_2(x,y)$. Then there is a finite set of 
numbers such that unless ${\beta \over \alpha}$ is in this set, $(-j_{\alpha,\beta}, p_{\alpha,\beta})
\leq \min((-j_1,p_1), (-j_2,p_2))$ under the lexicographic ordering.

One includes $\infty$ as a possible value of ${\beta \over \alpha}$ in Theorem 1.3. Also, note
that one can make appropriate generalizations of Theorem 1.3 for several functions.

We now give an idea of how Theorems 1.1 and 1.2 are proved. First consider the simple case where $D$ is a 
small disk centered at the origin and $S(x,y)$ and $f(x,y)$ are monomials $a_1x^{\alpha_1}y^{\beta_1}$ and 
$a_2x^{\alpha_2}y^{\beta_2}$. Then using an elementary
argument, one can evaluate $M_{S + tf, D}(\epsilon)$ directly to show Theorems 1.1 and 1.2. Already one might have to
exclude one value of $t$; this occurs if $\alpha_1 = \alpha_2$, $\beta_1 = \beta_2$, and $t = -{a_1 \over 
a_2}$. 
Next, suppose that instead of being monomials, $S(x,y)$ and $f(x,y)$ are comparable to monomials. That is,
suppose there are smooth functions $a_1(x,y)$ and $a_2(x,y)$, both nonvanishing at the origin, such that
$S(x,y) = a_1(x,y)x^{\alpha_1}y^{\beta_1}$ and $f(x,y) = a_2(x,y)x^{\alpha_2}y^{\beta_2}$. Then roughly
speaking one has the same behavior as in the monomial case. There is an added difficulty if $\alpha_1 = 
\alpha_2$ and $\beta_1 = \beta_2$, for in this case one must also have to consider the zeroes of 
$a_1(x,y) + ta_2(x,y)$. 

More generally, the strong form of resolution of singularities says that in the real-analytic case there 
is a coordinate change $\phi$ such that $S \circ \phi(x,y)$ and $f \circ \phi(x,y)$ are locally comparable 
to monomials in the above sense. However, the best one can automatically say about the Jacobian of this 
coordinate change is that it too is comparable to a monomial. Hence when looking at integrals one cannot 
automatically reduce to the above situations in general. 
Fortunately, in two dimensions there are substitutes for such resolution of singularities algorithms that
reduce to situations similar to where $S \circ \phi$ and $f \circ \phi$ are locally comparable 
to monomials, and which have determinant one. Even better, these algorithms hold for the general smooth case.

Specifically, in section 3, we will take a small disk $D$ centered at the origin and write
$D = \cup_{i=1}^n D_i$. On each $D_i$ there will be a coordinate change $\phi_i$ such that on each
$\phi_i^{-1} D_i$, the function $S \circ \phi_i$ is comparable to a monomial in a certain sense. This can be
done in such a way that $f \circ \phi_i$ is also comparable to a monomial, although we won't explicitly 
use this fact for certain technical reasons. Each $\phi_i(x,y)$ is of the form $(\pm x, \pm y - g_i(x))$,
and the domain $\phi_i^{-1} D_i$ is a "curved triangle" consisting of the points in $\phi_i^{-1}
(D_i)$ between two curves $y = p_i(x)$ and $y = q_i(x)$ such that $p_i(0) = q_i(0) = 0$ and $p_i(x^N)$ and
$q_i(x^N)$ are smooth for some $N$. Since each $\phi_i$ now has Jacobian determinant $\pm 1$, in examining
$M_{S,D}(\epsilon)$ one can switch to 
considering $S \circ \phi_i$ and $f \circ \phi_i$ on the set $\phi_i^{-1} D_i$. Although these two functions
aren't strictly speaking comparable to monomials, there are enough similarities with that situation such that 
after some effort one can prove Theorem 1.2. One has to exclude finitely many values of $t$ for each 
$D_i$ to avoid cancellations such as in the monomial case. 

The idea of dividing into curved triangles related to the singularities of $S(x,y)$ to simplify the behavior
of integrals related to $S(x,y)$ goes a while back. It was 
used in the various Phong-Stein papers on oscillatory integral operators such as [PS] and then in the 
author's earlier work such as [G2]-[G3].
The Phong-Stein papers use curved triangles deriving from Puiseux expansions of real-analytic 
functions, while [G2]-[G3] uses explicit resolution of singularities algorithms such as in this paper.
Since the problems being considered here are rather different from the earlier problems, we will derive
from first principles a resolution of singularities theorem amenable to the situations at hand. 

Proving the stronger result of Theorem 1.1 requires additional ideas. In fact, if our goal was only to 
prove Theorem 1.2 and its consequences, then section 4 would be noticeably shorter. To get the sharper
estimates of Theorem 1.1, we will draw on the results of [G1]. Specifically, we first put $S(x,y)$ into 
what in [G1] are called superadapted coordinates. These are a generalization of the notion of adapted 
coordinates of [V]. Then 
we apply the resolution of singularities algorithm of section 3 to $S(x,y)$,
getting the resulting $D_i$. One next focuses on the $D_i$ which give the dominant terms in $M_{S,D}
(\epsilon)$. For these $D_i$, one subdivides further into sets $D_{ij}$. This will be a coarse 
subdivision related to the resolution of singularities of $S(x,y) + tf(x,y)$; however, one does not have to 
use the full resolution of singularities theorem here. For the $D_{ij}$ 
that give the largest contribution to $M_{S,D}(\epsilon)$, one uses estimates from section 2  
related to one-dimensional Van der Corput-type lemmas to prove the sharp estimates of Theorem 1.1.

On the remaining $D_{ij}$, as well as the $D_i$
that do not give dominant terms of $M_{S,D}(\epsilon)$, one now applies the full resolution of singularities
theorem to $f \circ \phi_i(x,y)$. The resulting functions, call them $f \circ \Phi_{ij}$ and $S \circ 
\Phi_{ij}$, are now comparable to monomials in the new coordinates, and the considerations used for
Theorem 1.2 can now be used. Because the contributions here are error terms for $M_{S,D}(\epsilon)$ (that is,
they give
higher powers of $\epsilon$ than the estimates sought), we do not have to worry about constants here
if $\epsilon$ is small enough, which we will see we can assume. However, we still have to exclude finitely
many values of $t$ as in Theorem 1.2; it is conceivable that for such $t$ the power of epsilon appearing in
such error terms 
becomes as small or smaller than the desired power for $M_{S,D}(\epsilon)$. Note that this phenomenon does
not appear in the real-analytic results [K1]-[K3] or [PSSt]. The author does not know if this is a result
of the weaker assumptions of Theorem 1.1, or if one can avoid excluding finitely many values of $t$ in the
context of Theorem 1.1 with an additional argument. As indicated above, this cannot be avoided in the context
of Theorem 1.2.

\noindent {\bf Oscillatory Integrals.}

Let $S(x,y)$ be a smooth function on a neighborhood $D$ of the origin with $S(0,0) = 0$. Suppose
$\phi(x,y)$ is a real-valued smooth function supported in $D$. We consider the oscillatory integral 
defined by
$$J_{S,\phi}(\lambda) = \int_{\R^2} e^{i \lambda S(x,y)} \phi(x,y)\,dx\,dy \eqno (1.5)$$
Here $\lambda$ is a real parameter and we are interested in the behavior of $J_{S,\phi}(\lambda)$ as
$|\lambda| \rightarrow \infty$. Since  $J_{S,\phi}(-\lambda)$ is the complex conjugate of $J_{S,\phi}
(\lambda)$, one only needs to consider the situation as $\lambda \rightarrow \infty$. As is well-known,
the analysis of $(1.5)$ is closely related to the analysis of the sublevel areas above. Specifically,
in the real-analytic case, if $S(0,0) = 0$ and $\nabla S(0,0) = 0$ then in analogy to $(1.2)$ we have 
nontrivial asymptotics of the form
$$J_{S,\phi}(\lambda)= D_{\phi} \lambda^{-j}(\ln \lambda)^p + o(\lambda^{-j}(\ln \lambda)^j) \eqno (1.6)$$
Here we write $(1.6)$ such that $(-j,p)$ is maximal (under the lexicographic ordering) such that $D_{\phi}$ is nonzero for at least one
$\phi$ in any sufficiently small neighborhood of the origin. We refer to $j$ as the oscillatory index of
$J_{S,\phi}$ and $p$ as its multiplicity. Using well-known arguments from [AGV], it can be shown that
(in the real-analytic case) the oscillatory index is equal to the growth index, with the multiplicity 
the same in both cases, unless
there is a coordinate system near the origin in which $S(x,y) = x^2 - y^2$. In this case, the growth 
and oscillatory indices are both 1, with the multiplicity of the growth index being 1 and the multiplicity
of the oscillatory index being zero. 

Furthermore, Karpushkin's methods work for the oscillatory integral
case as well, and in [K1]-[K2] analogues for oscillatory integrals $(1.6)$ to his above-mentioned theorem
on sublevel areas are proven. 

Using some  results of [G1], Theorem 1.3 directly implies analogues for oscillatory integrals. To see why 
this is the case, we first give some background from [G1] which will also be used in the proof of 
Theorem 1.1. Suppose $S(x,y)$ is a smooth function on a neighborhood of the origin such that
$S(0,0) = 0$ and $\nabla S(x,y)$. Write the Taylor series of $S(x,y)$ as $\sum_{ab} s_{ab} x^ay^b$ and 
denote the Newton distance of $S(x,y)$ by $d$. Then 
it is proven in [G1] that there is a smooth coordinate change taking the origin to itself, such that
after the coordinate change $S(x,y)$ is in "superadapted coordinates", which means the following. 

\noindent {\bf Definition.} For a compact edge $e$ of $N(S)$ with equation $a + mb = \alpha$, let
$S_e(x,y)$ denote the sum of all terms of $\sum_{a,b} s_{ab} x^ay^b$ with $a + mb = \alpha$. $S(x,y)$ is 
said to be in {\it superadapted coordinates} if for any compact edge $e$ of $N(S)$ intersecting the 
bisectrix, any zero of $S_e(1,y)$ or $S_e(-1,y)$ has order less than $d(S)$.

In [G1] it is proven that in superadapted coordinates, the growth index $j$ is given by ${1 \over d}$. 
The multiplicity $p$ is 1 if and only if the bisectrix intersects $N(S)$ at a vertex. It is also shown that
$d = 1$ if $S(x,y)$ is Morse and $d > 1$ otherwise. It is further shown that in superadapted coordinates,
on a small enough neighborhood $U$ of the origin for these values of $j$ and $p$ one has $J_{S,\phi}(\lambda) < C 
\lambda^{-j}(\ln \lambda)^p $ for any $\phi$ supported in $U$, and that for any $(-j',p') < (-j,p)$ 
(with respect to the lexicographic ordering) there is some $\phi$ supported 
on $U$ for which the estimate $J_{S,\phi}(\lambda) < C \lambda^{-j'}(\ln \lambda)^{p'}$ does not 
hold for {\it any} $C$. 

Thus in the non-Morse case, for a general smooth function it makes sense to define the
oscillatory index and multiplicity to be the same as the growth index and multiplicity. In the Morse
case one defines them to be inherited from the Morse coordinates. Note that this definition agrees with
the old definition for the real-analytic case. Note also that the analogue of Theorem 1.3 for 
oscillatory integrals follows immediately from Theorem 1.3. Since
the two types of Morse critical points have the same oscillatory indices and multiplicities, the case 
where $S_1(x,y)$ and $S_2(x,y)$ both have Morse critical points at the origin may be included in the 
oscillatory integral result: 

\noindent {\bf Theorem 1.4.} Suppose $S_1(x,y)$ and $S_2(x,y)$ are smooth functions on a 
neighborhood of the origin with critical points at the origin. Let $(j_1,p_1)$ and $(j_2,p_2)$ be their 
oscillatory indices and multiplicities at the origin. Let $j_{\alpha,\beta}$ and $p_{\alpha,\beta}$ be the 
oscillatory index and multiplicity of $\alpha S_1(x,y) + \beta S_2(x,y)$. Then there is a finite set of real
numbers such that unless ${\beta \over \alpha}$ is in this set, then $(-j_{\alpha,\beta}, 
p_{\alpha,\beta}) \leq \min((-j_1,p_1), (-j_2,p_2))$ under the lexicographic ordering.

Again, here we include $\infty$ as a possible value of ${\beta \over \alpha}$. Note that Theorem 1.4 pertains
to integrals of the form $\int_{R^2} e^{i\alpha S_1(x,y) + i\beta S_2(x,y)}\phi(x,y)\,dx\,dy$, which are
tied to Fourier transforms of surface-supported measures. 

For oscillatory integrals, getting constants depending on $S(x,y)$ and not $f(x,y)$ in analogy with
Theorem 1.1 is difficult for a few reasons. For one, since the integrands of oscillatory integrals have
both positive and negative values, even if one had precise one-dimensional Van der Corput lemmas for
the oscillatory integral case, averaging the resulting estimates over a 
second dimension might give extra cancellation that needs to be taken into account. Secondly, and perhaps
more importantly, applying one-dimensional Van der Corput lemmas on the integrands of $J_{S,\phi}$ will
result in bounds depending on the supremum of
$|\nabla\phi|$, and such upper bounds do not necessarily behave well under coordinate changes. However, 
if the perturbed function $S + tf$ is real-analytic there are explicit formulas from [AGV] for 
transforming the estimates of Theorem 1.1 into explicit estimates for large enough $\lambda$, and we get
the following.

\noindent {\bf Theorem 1.5.} Suppose $S(x,y)$, $f(x,y)$, $l$, $\delta$, $j$, $p$, and $D$ are as in
Theorem 1.1. There is a $B_S > 0$ such
that for all but finitely many $|t| < 1$, if $\phi$ is supported in $D$ and
$S + tf$ is real-analytic, then for sufficiently large $\lambda$ we have  
$$|J_{S + tf,\phi}(\lambda)| <  B_S||\phi||_{\infty}\lambda^{-j}(\ln \lambda)^p\eqno (1.7)$$

\noindent {\bf Proof.} We may assume $S(x,y)$ is not Morse since the Morse case follows from the
explicit asymptotic expansions known in the Morse situation. By Theorem 1.2 of [G1], the initial 
coefficient $D_{\phi}$ of the asymptotics $(1.6)$ for $S + tf$ satisfies 
$$|D_{\phi}| \leq j\Gamma(j) \lim_{\epsilon \rightarrow 0}\bigg|{I_{S + tf,\phi}(\epsilon) 
\over \epsilon^{j}(\ln \epsilon)^p}\bigg| \eqno (1.8)$$
Here $I_{S + tf,\phi}(\epsilon) = \int_{|S + tf| < \epsilon} \phi(x,y)\,dx\,dy$. Note that we have
$$\bigg|{I_{S + tf,\phi}(\epsilon) \over \epsilon^{j}(\ln \epsilon)^p}\bigg| \leq ||\phi||_{\infty}
\bigg|{M_{S + tf,D}(\epsilon) \over \epsilon^{j}(\ln \epsilon)^p}\bigg|$$
Since Theorem 1.1 holds for $S + tf$ on $D$, the above is at most
$$C_S ||\phi||_{\infty}$$
Thus since $j\Gamma(j) \leq 2$, the limit in $(1.8)$ is at most $2C_S ||\phi||_{\infty}$.
By taking $\lambda$ sufficiently large that the other terms of the asymptotics
$(1.6)$ are small in comparison, we obtain $(1.7)$ with $B_S = 3C_S$ and the theorem follows. 

\noindent {\bf 2. Lemmas about curved triangles and one-dimensional Van der Corput Lemmas}

\noindent We make extensive use of the classical Van der Corput lemma throughout this paper. Although 
we don't need very sharp constants, to simplify our notation we use the following version that follows 
from [R]. We refer to [CaCWr] for more information on this general subject.

\noindent {\bf Lemma 2.1.} ([R]) Suppose for a positive integer $k$, $f(t)$ is a $C^k$ function on an interval 
$I$ such that for some positive constant $c$, $|{d^k f \over dt^k}| > c k!$ on $I$. Then for any 
$\epsilon > 0$ we have
$$|\{t \in I: |f(t)| < \epsilon\}| \leq \min(|I|,  4 c^{-{1 \over k}} \epsilon^{1 \over k}) 
\eqno (2.1)$$
An immediate consequence of Lemma 2.1 that will be useful in analyzing our functions on curved triangles
is the following.

\noindent {\bf Lemma 2.2.} Let $A_{m,N}$ denote the set $\{(x,y): 0 < x < x_0,\,\,0 < y < Nx^m\}$,
where $m > 0$. Suppose $g(x,y)$ is a function on $A_{m,N}$ such that $|\partial_y^{\beta}
g(x,y)| > a\beta ! x^{\alpha}$ on $A_{m,N}$, where $a,\alpha > 0$ and $\beta$ is a positive integer. 
Then for a fixed $x$ we have
$$|\{y: (x,y) \in A_{m,N},\,\,|g(x,y)| < \epsilon\}| \leq 4 \{y: (x,y) \in A_{m,N},\,\, 
ax^{\alpha}y^{\beta} < \epsilon\}| \eqno (2.2)$$
Lemma 2.2 will often be used in conjunction with the following lemma.

\noindent {\bf Lemma 2.3.} Let $h(x,y) = ax^{\alpha}y^{\beta}$ for some $a,\alpha, \beta \geq 0$. Let
$A_{m,N}$ be as in Lemma 2.2. 

\noindent {\bf a)} If $\beta > \alpha$, then there are constants $c,C > 0$ depending on $a,m,\alpha$, and 
$\beta$ such that for sufficiently small $\epsilon$ we have 
$$c x_0^{\beta - \alpha \over \beta}\epsilon^{1 \over \beta} < |\{(x,y) \in A_{m,1}: |h(x,y)| 
< \epsilon\}| < C x_0^{\beta - \alpha \over \beta}\epsilon^{1 \over \beta} \eqno (2.3)$$
\noindent {\bf b)} If $\beta = \alpha$, then there are constants $c,C > 0$ depending on $a, \alpha, m$,
and $\beta$ such that for sufficiently small $\epsilon$ we have the estimate 
$$c\epsilon^{1 \over \beta} |\ln \epsilon| < 
|\{(x,y) \in A_{m,1}: |h(x,y)| < \epsilon\}| <C \epsilon^{1 \over \beta} |\ln 
\epsilon| \eqno (2.4)$$
\noindent {\bf c)} If $\beta < \alpha$, then there are constants $c,C > 0$ depending on $a,m,\alpha$, and 
$\beta$ such that for sufficiently small $\epsilon$ we have 
$$c N^{\alpha - \beta \over \alpha + m\beta}\epsilon^{{m + 1 \over \alpha + m\beta}}< |\{(x,y) 
\in A_{m,N}: |h(x,y)| < \epsilon\}| < C N^{\alpha - \beta \over \alpha + m\beta}\epsilon^{{m + 1 \over 
\alpha + m\beta}} \eqno (2.5)$$
\noindent {\bf Proof:} Viewing $|\{(x,y) \in A_{m,N}: |h(x,y)| <\epsilon\}|$ as the integral of the 
characteristic function of $\{|h(x,y)| <\epsilon\}$ over $A_{m,N}$, we change variables twice, first
by replacing $y$ by $x^my'$ and then by replacing $x$ by $(x')^{1 \over m+1}$. We obtain that $|\{(x,y) 
\in A_{m,N}: |h(x,y)| <\epsilon\}|$ is given by 
$$(m+1)^{-1}|\{(x',y') \in [0,x_0^{m+1}] \times 
[0,N]: a(x')^{\alpha + m\beta \over m+1}(y')^{\beta} < \epsilon\}| \eqno (2.6)$$
The $x'$-measure in $(2.6)$ for fixed $y'$ is given by $\min(x_0^{ m+1}, c\epsilon^{m + 1 \over
\alpha + m\beta}
(y')^{-{m\beta + \beta \over \alpha + m \beta}})$. Note that the two terms in the minimum are equal at
$y = y_0' = c'x_0^{-{\alpha + m\beta \over \beta}}\epsilon^{1 \over \beta}$. Also note that the power
${m\beta + \beta \over \alpha + m \beta}$ of $y'$ is greater than 1 if $\beta > \alpha$, and
less than 1 if $\beta < \alpha$. In the former case, if $\epsilon$ is sufficiently small 
the measure of the set $(2.6)$ is comparable to
the portion where $y' < y_0'$, given by $y_0'x_0^{ m+1}$ = $c'x_0^{\beta - \alpha \over 
\beta}\epsilon^{1 \over \beta}$ and thus the formula $(2.3)$ holds. If $\beta = 
\alpha$, then the exponent is exactly 1, and one obtains the additional logarithmic factor of $(2.4)$.
Lastly, if $\beta < \alpha$, the measure of the set $(2.6)$ is comparable to the measure of the part 
where $N > y' > {N \over 2}$, giving $(2.5)$ for small enough $\epsilon$. This completes the proof of 
the lemma.

\noindent We will make frequent use of the next lemma in conjunction with the above lemmas. 

\noindent {\bf Lemma 2.4.} Suppose $R(x,y)$ is a smooth function on a neighborhood of $(0,0)$ such that
$N(R)$ has a vertex at $(c,d)$. Let $r_{cd}x^cy^d$ be the associated term of the Taylor expansion of
$R(x,y)$ at $(0,0)$. Then the following hold.

\noindent {\bf a)} If $(c,d)$ is the intersection of two compact edges of $N(R)$, with equations $a + mb
= \alpha$ and $a + m'b = \alpha'$ respectively with $m' > m$, then for any $\delta > 0$ there is an
$r > 0$ and a $\xi > 0$ such that for $0 < x < r$ and $\xi^{-1}x^{m'} < y < \xi x^m$ one has 
$$|R(x,y) - r_{cd}x^cy^d| < \delta|r_{cd}|x^cy^d \eqno (2.7)$$
\noindent {\bf b)} If $(c,d)$ is the intersection of the horizontal ray of $N(R)$ with a compact edge 
with equation $a + mb = \alpha$, then for sufficiently small $\eta > 0$, for any $\delta > 0$ there is an 
$r > 0$ and a $\xi > 0$ such that for $0 < x < r$ and $x^{1 \over \eta} < y < \xi x^m$ equation $(2.7)$ 
holds.

\noindent {\bf c)} In the setting of case b), if we also have that $d = 0$ then for any $\delta > 0$ 
there is an $r > 0$ and a $\xi > 0$ such that for $0 < x < r$ and $0 < y < \xi x^m$ equation $(2.7)$ 
holds.

\noindent {\bf Proof.} We start with a). Let the Taylor expansion of $R(x,y)$ at the origin be written as 
$\sum_{ab} r_{ab}x^ay^b$. For a large $M$ we may write 
$$R(x,y) -  r_{cd}x^cy^d = \sum_{c \leq a < M,\,\,d \leq b < M,\,\,(a,b) \neq (c,d)}r_{cd}x^ay^b$$
$$ + \sum_{a < c,\,\,d < b < M ,\,\,a + mb \geq \alpha}r_{ab}x^ay^b + 
\sum_{c < a < M,\,\,b < d,\,\,a + m'b \geq \alpha '}r_{ab}x^ay^b + E_M(x,y) \eqno (2.8)$$
Here $E_M(x,y)$ satisfies 
$$|E_M(x,y)| < C(|x|^M + |y|^M) \eqno (2.9)$$
The first sum in $(2.7)$ can be made less than ${\delta \over 4}|r_{cd}|x^cy^d$ in absolute value by shrinking
the radius of $D$ appropriately. As for the second sum, if one changes coordinates from $(x,y)$ to 
$(x,y')$, where $y' = x^{m}y$, then $(x,y') \in [0,1] \times [0,\xi]$.
Observe that a given term $r_{ab}x^ay^b$ of the second sum becomes $r_{ab}x^{a + mb}(y')^b$. Since
$a + mb \geq \alpha$ and $b > d$ in each term in the second sum, the entire sum can be written as
$x^{\alpha}(y')^d(y'f(x,y'))$ for some $f(x,y')$ which is a polynomial in $y'$ and a fractional power of
$x$. Thus by shrinking $\xi$ appropriately, since $y' < \xi$ the sum can be made of absolute value 
less than ${\delta \over 4}|r_{cd}|
x^{\alpha}(y')^d$. Note that $r_{cd}x^cy^d = r_{cd}x^{c + dm}(y')^d$, and this is equal to 
$r_{cd}x^{\alpha}(y')^d$ since $(c,d)$ is on the edge with equation $a + mb = \alpha$. Thus
by choosing $\xi$ sufficiently small, we can have that the second sum is of absolute value at most 
${\delta \over 4}|r_{cd}|x^cy^d$ (in the original coordinates.) These are the bounds we need.

The third sum is dealt with in exactly the same way, reversing the roles of the $x$ and $y$ axes. Lastly,
since $x^m > y > x^{m'}$ the error term $E_M(x,y)$ can be made less than ${\delta \over 4} |r_{cd}|x^cy^d$ in 
absolute value for small $x$ if $M$ is chosen sufficiently large. Putting these all together,
we obtain $|R(x,y) - r_{cd}x^cy^d| < \delta|r_{cd}|x^cy^d$ as needed. This completes the proof of part a).
of the lemma.

Moving on now to part b), we once again look at the expression $(2.8)$. This time there is no third 
sum, and all the other terms can be bounded exactly as in part a); the condition that $y > x^{1 \over 
\eta}$ for some small $\eta$ ensures that for large enough $M$ the error term will still be smaller
than ${\delta \over 4}|r_{cd}|x^cy^d$ using the inequality $(2.9)$. Moving now to c), we again examine
$(2.8)$. Once again there is no third sum, and the first two terms are bounded as they were in part a).
This time the error term is bounded by $C(|x|^M + |y|^M) \leq C'(|x|^M + |x|^{Mm})$. Since $d = 0$, 
by taking $M$ large enough the error term can be made less than ${\delta \over 4}|r_{cd}|x^cy^d$ 
$= {\delta \over 4}|r_{cd}|x^c$ for small enough $x$. This gives us part c) and completes the proof of
Lemma 2.4. 

\noindent {\bf 3. The resolution of singularities algorithm.}

In this section we prove the version of two-dimensional resolution of singularities we need for the
analysis in section 4. In keeping with the philosophy of [G2] as well as its antecedents such as
[PS] or [V], it involves dividing a 
neighborhood of the origin into "curved triangles" each of which has a coordinate system in which 
$S(x,y)$ behaves like a monomial in an appropriate sense. The theorem we use is the following.

\noindent {\bf Theorem 3.1.} Suppose $S(x,y)$ is a smooth function defined on a 
neighborhood of the origin with $S(0,0) = 0$ such that the Taylor expansion $\sum_{ab} s_{ab}x^ay^b$
of $S(x,y)$ at the origin has at least one nonvanishing term. Then for any sufficiently small $\eta > 0$,
and any sufficiently small disk $D$ centered at the origin, we may, up to a set
of measure zero, write $D \cap \{(x,y): |y| < |x|^{\eta}\}$  as a finite union $\cup_i D_i$, where 
the $D_i$ have the following properties. 

\noindent Let $M$ denote the difference between the
$y$ coordinates of the uppermost and lowermost vertices of $N(S)$. If $N(S)$ has one vertex let
$M = 1$. Then there is a positive integer $M'$ depending on $M$ such that for each $i$ there is a function
$g_i(x)$ with $g_i(x^{M'})$ smooth and a function $\phi_i(x,y)$ of the form $(\pm x, \pm y - g_i(x))$
such that 
$f \circ \phi_i(x,y)$ is approximately a nonzero monomial $b_i x^{\alpha_i}y^{\beta_i}$ ($\alpha_i, 
\beta_i \geq 0, \beta_i$ an integer) on the curved triangle $D_i' = \phi_i^{-1}D_i$ in the following 
sense. 

\noindent {\bf a)} The set $D_i'$ is of the form $\{(x,y) \in \phi_i^{-1}(D): x > 0,\,\,
f_i(x) < y < F_i(x) \}$ where each $f_i(x^{M'})$ and $F_i(x^{M'})$ are smooth. The initial term of the
Taylor expansion of $F_i(x)$ is of the form $A_i x^{N_i}$, where $A_i, N_i > 0$. The function $f_i(x)$ 
is either the zero function or has a Taylor series with initial term $a_i x^{n_i}$ where $a_i, n_i
> 0$ and $n_i > N_i$.

\noindent {\bf b)} If $\beta_i = 0$, then $f_i(x)$ is the zero function and there are positive constants 
$c_i$ and $C_i$ such that on $D_i'$ one has the estimates
$$c_ix^{\alpha_i} < S \circ \phi_i(x,y) < C_ix^{\alpha_i}$$
\noindent {\bf c)} If $\beta_i > 0$, then we can write $S = S_1 + S_2$ as follows. $S_2$ has a zero of 
infinite order at $(0,0)$ and is the zero function if $S$ is real-analytic. Also, $S_2(x^{M'},y)$ is a 
smooth function. As for $S_1$, for any preselected $\delta > 0$ we can arrange that the 
decomposition is such that for some nonzero constant $b_i$, for any $0 \leq k \leq \alpha_i$, 
$0 \leq l \leq \beta_i$ on $D_i'$ we have
$$|\partial_x^k\partial_y^l S \circ \phi_i(x,y) - b_i\alpha_i(\alpha_i - 1)...(\alpha_i - k + 1)
\beta_i(\beta_i - 1)....(\beta_i - l + 1)x^{\alpha_i - k}y^{\beta_i - l}| $$
$$\leq \delta|b_i|x^{\alpha_i - l} y^{\beta_i - k} \eqno (3.1)$$
\noindent {\bf d)} The total number of sets $D_i$ can be bounded in terms of $M$.

\noindent The following corollary will follow immediately from the proof of Theorem 3.1.

\noindent {\bf Corollary 3.2.} If $S(x^p,y)$ is a smooth function for some positive
integer $p$, then Theorem 3.1 still holds, except the exponent $M'$ is replaced by $pM'$. 

\noindent {\bf Proof of Theorem 3.1.} We first dispense with the case where $N(S)$ has exactly one vertex $(a,b)$. Let
$s_{ab}x^ay^b$ be the corresponding term of the Taylor expansion of $S(x,y)$. Thus for any $N$ we have 
$S(x,y) = s_{ab}x^ay^b
+ O(|x|^N + |y|^N)$. We divide a small disk $D$ centered at the origin into 8 curved triangles by 
slicing along the $x$ and $y$ axes as well as along the lines $y = \pm x$. We make these triangles the 
$D_i$'s with each $\phi_i(x) = (\pm x, \pm y)$ or $(\pm y, \pm x)$ and each $f_i(x) = 0$. Then condition 
c) above automatically holds and we are done.

Thus we now assume that $N(S)$ has multiple vertices. Hence $N(S)$ has at least one compact edge. We 
write the equations of these edges as $a + m_jb = \alpha_j$, where $m_1 > m_2 > ... > m_n$. Clearly it
suffices to divide the $x > 0, y > 0$ portion of a small disk $D$ centered at the origin according to the
lemma, so we restrict our attention to $D_0 = D \cap \{x > 0, y > 0\}$. For a small $\xi > 0$ to be 
determined, up to a set of measure zero we write 
$D_0 \cap \{y < x^{\eta}\}$ as the finite union of sets $U_j$ and $T_j$, where the $U_j$ are all
possible sets of the form $\{(x,y) \in D_0: \xi x^{m_j} < y < \xi^{-i} x^{m_j}\}$ and the $T_j$ are all
possible sets of the form $\{(x,y) \in D_0: \xi^{-i} x^{m_{j+1}} < y < \xi x^{m_j}\}$ as well as the
sets $T_n = \{(x,y) \in D_0: 0 < y < \xi x^{m_n}\}$ and $T_0 = \{(x,y) \in D_0: \xi^{-i} x^{m_1} < y
< x^{\eta}\}$. (We assume $\eta < m_1$). Observe that each $T_j$ corresponds to a unique vertex of 
$N(S)$. We will turn each $T_j$ into
one of the $D_i$'s with the associated $\phi_i$ just the identity map. As for the $U_j$, we will
subdivide each $U_j$ further into $V_{jk}$ and $U_{jl}$, where the $U_{jl}$ will also become 
$D_i$'s for which part b) of the theorem holds, and where the $V_{jk}$ will undergo further subdivisions and
coordinate changes. 

We start with the $T_j$'s for $j \neq n$. Let $(c,d)$ denote the vertex of $N(S)$ corresponding to
$T_j$ and let $s_{cd}x^cy^d$ denote the associated term of the Taylor expansion. 
Then by Lemma 2.4 applied to $S(x,y)$ and its various $y$ partials, for any $\delta > 0$ we can choose 
$\xi$ such that if $D$ is sufficiently small then on $T_j$, we have the inequality 
$|S(x,y) - s_{cd}x^cy^d| < \delta|s_{cd}|x^cy^d$ as well as its analogues for any $\partial_x^k\partial
_y^lS(x,y)$ for $k \leq c$, $l \leq d$. Thus a) and c) of the theorem hold with $\phi_i$ the identity map, 
which is what we need for these $T_j$.

Next, we look at $T_n = \{(x,y) \in D_0: 0 < y < \xi x^{m_n}\}$. This time we cannot apply Lemma 2.4 
automatically. By a famous theorem of Borel
(see [H] for a proof), one can let $s_0(x,y)$ be a smooth function on a neighborhood of the origin whose 
Taylor expansion at the origin is given by $\sum_{ab} s_{ab}x^ay^{b-d}$. Then Lemma 2.4c) applies to
$s_0(x,y)$ and we can assume that on $T_n$ we have 
$$|s_0(x,y) - s_{cd}x^c| < \delta|s_{cd}x^c| \eqno (3.2)$$
Let $s_1(x,y) = y^d s_0(x,y)$. Then $s_1(x,y)$ has the same Taylor expansion at the origin as $S(x,y)$
and is equal to $S(x,y)$ when $S(x,y)$ is real-analytic. We also have the desired inequality
$$|s_1(x,y) - s_{cd}x^cy^d| < \delta|s_{cd}x^cy^d| \eqno (3.3)$$
Note also that the analogues of $(3.3)$ for the $x$ and $y$ partials also hold; for example, the Newton
polygon of $s_0(x,y)$
is such that taking a $y$ derivative of $s_0(x,y)$ incurs a factor of at most $Cx^{-m_n}$ which is
much smaller than ${1 \over y}$ on $T_n$ if $\xi$ is appropriately small. Thus we have a) and c) of 
Theorem 3.1 and we are done with the analysis of the $T_j$'s.

Next, we proceed to the analysis of the $U_j$'s. Let $S_{m_j}(x,y)$ denote the sum of the terms $s_{ab}x^ay^b$
of the Taylor expansion lying on the edge with equation $a + m_jb = \alpha_j$. Note that $S_{m_j}(x,y)$
is a polynomial and is the sum of $s_{ab}x^ay^b$ minimizing $a + m_jb$. Let $r_{j1} < ... < r_{jN_j}$ 
denote the positive zeroes of $S_{m_j}(1,y)$ if there are any. Define $V_{jk} = \{(x,y) \in U_j: 
(r_{jk} - \xi)x^{m_j} < y < (r_{jk} + \xi)x^{m_j}\}$. As long as $\xi$ is small enough, we may
write $U_j - \cup_k V_{jk}$ as a union $\cup_l U_{jl}$ where each $U_{jl}$ is of the form $\{(x,y) \in
U_j: (r_{j\,k+1} - \xi)x^{m_j} < y < (r_{jk} - \xi)x^{m_j}\}$, $\{(x,y) \in U_j: \xi x^{m_j} < y < (r_1 - \xi)x^{m_j}\}$, or 
$\{(x,y) \in U_j: (r_{N_j} + \xi)x^{m_j} < y < \xi^{-1}x^{m_j}\}$. In the case that $S_{m_j}(1,y)$ has
no positive zeroes, there are no $V_{kl}$'s and we just set $U_{j1} = U_j$.

On each $U_{jl}$, $S(x,y)$ is already in the form required in Theorem 3.1. To see this, one does the 
coordinate change $(x,y) = (x,x^{m_j}y')$, turning $U_{jl}$ into a set $U_{jl}'$ contained in 
$[0,\eta] \times [a_{jl},b_{jl}]$, where $\eta$ denotes the radius of the disk $D$ and where 
$a_{jl}, b_{jl} > 0$. Note that $S_{m_j}(1,y')$ has no zeroes on $[a_{jl},b_{jl}]$. In the new coordinates,
the finite Taylor expansion $S(x,y) = \sum_{a,b < M}s_{ab}x^ay^b + O(|x|^M + |y|^M)$ becomes of the form
$$S(x,x^{m_j}y') = x^{\alpha_j}S_{m_j}(1,y') + x^{\alpha_j + \zeta_j}f(x,y') + O(|x|^M + |x|^{mM}|y'|^M) \eqno
(3.4)$$
Here $f(x,y')$ is a polynomial in $y'$ and a fractional power of $x$ and $\zeta_j > 0$. Since $S_{m_j}
(1,y')$ has no zeroes on $[a_{jl},b_{jl}]$, there are $C_{jl}, \epsilon_{jl} > 0$ such that 
$C_{jl}x^{\alpha_j} > |x^{\alpha_j}S_{m_j}(1,y')| > 
\epsilon_{jl}x^{\alpha_j}$ on $U_{jl}'$. Furthermore, if $\eta$ is sufficiently small and $M$ is 
sufficiently large, we have that $|x^{\alpha_j + \zeta}f(x,y')|$ and the $O(|x|^M + |x|^{mM}|y'|^M)$ terms
are both less than $\min ({C_{jl} \over 4}x^{\alpha_j},{\epsilon_{jl} \over 4}x^{\alpha_j})$. As a result, 
shrinking the disk $D$ if necessary we can assume that on $U_{jl}'$ we have
$${C_{jl} \over 2}x^{\alpha_j} > |S(x,x^{m_j}y')| > {\epsilon_{jl} \over 2}x^{\alpha_j} \eqno (3.5)$$
Translating back into the coordinates of $U_{jl}$, this means that on $U_{jl}$ we have
$${C_{jl} \over 2}x^{\alpha_j}> |S(x,y)| > {\epsilon_{jl} \over 2}x^{\alpha_j} \eqno (3.6)$$
Thus on $U_{jl}$, $S(x,y)$ satisfies the conditions of Theorem 3.1b) with $\beta_i = 0$, if we let the
coordinate change $\phi_i$ associated with $U_{jl}$ be given by $(x,y) \rightarrow (x,y + c_{jl}x^m)$,
where $y = c_{jl}x^{m_j}$ denotes the equation of the lower boundary curve of $U_{jl}$. This completes 
the analysis of the $U_{jl}$.

We now move to the analysis of any $V_{jk}$ that may exist. Let $o_{jk}$ denote the order of the zero of
$S_{m_j}(1,y)$ at $y = r_{jk}$. The idea is as follows. Using ideas from resolution of singularities, we will
do a coordinate change of the form $(x,y) \rightarrow (x,y - r_{jk}x^{m_j} +$ higher order terms$)$ such
that in the new coordinates $S(x,y)$ becomes a function whose analogues to the zeroes $r_{jk}$ each has
order $< o_{jk}$. Thus after at most $\max_{j,k} o_{jk}$ iterations, there will no longer be any sets
$V_{jk}$ and we will have divided $D_0$ into sets each of which is a $U_{jl}$ or $T_j$ in some 
iteration. Since each coordinate change will be of the form $(x,y) \rightarrow (x, \pm y - g(x))$, the
composition of finitely many such coordinate changes is of the desired form. (We can get $-y$ as well as
$y$ since after one of these coordinate changes the resulting set lies in both the upper right and lower
right quadrants). By the above analysis of $S(x,y)$ on the $T_j$ and $U_{jl}$, the resulting $S \circ
\phi_i(x,y)$ will satisfy Theorem 3.1 as needed. 

The coordinate change on $V_{jk}$ is chosen in the following way. We once again switch to the $(x,y')$
coordinates and make use of $(3.4)$. Let $V_{jk}'$ denote $V_{jk}$ in the new coordinates and define 
$s_j(x,y') = {S(x,x^{m_j}y') \over x^{\alpha_j}}$, where recall $\alpha_j = a + m_jb$ for $(a,b)$ on
$e_j$. We claim that $s_j(x^N,y')$ is a smooth function of $x$ and $y'$ on $V_{jk}'$ for some
$N$. To see this, observe that by $(3.4)$, for any $(p,q)$ we have
$$\partial_x^p\partial_{y'}^q(s_j(x^N,y')) = \partial_x^p\partial_{y'}^q (S_{m_j}(1,y')) + 
\partial_x^p\partial_{y'}^q(x^{N\zeta_j}f(x,y')) $$
$$+ O(|x|^{N(M-p)} + |x|^{N(mM - p)}|y'|^{M - q}) \eqno(3.7)$$
The magnitude of the error term here follows from corresponding bounds on the derivatives of the error
term in the partial Taylor expansion $S(x,y) = \sum_{a,b \leq M} s_{ab}x^ay^b + O(|x|^M + |y|^M)$. Thus
as long as $N\zeta$ is an integer, 
by taking $M$ large enough we have that $\partial_x^p\partial_{y'}^q s_j(x^N,y')$ exists and is continuous
on $V_{jk}'$. Hence $s_j(x^N,y)$ is a smooth function of $x$ and $y$ as needed.

Furthermore, $\partial_{y'}^{o_{jk}}s_j(x,r_{jk}) \neq 0$, while $\partial_{y'}^{o_{jk}-1}s_j(x,r_{jk}) \neq 0$.
Hence the implicit function theorem (applied to $s_j(x^{N},y')$) says that there is a function 
$t_{jk}(x)$ for small $x$ with $t_{jk}(x^N)$ smooth such that $t_{jk}(0) = r_{jk}$ and
$$\partial_{y'}^{o_{jk}-1}s_j(x,t_{jk}(x)) = 0 \eqno (3.8)$$
More generally, we also have 
$$c|y' - t_{jk}(x)| \leq |\partial_{y'}^{o_{jk}-1}s_j(x,y')| \leq C|y' - t_{jk}(x)| \eqno (3.9)$$
Translating in terms of $S(x,y)$, on $V_{jk}$ we have 
$$\partial_y^{o_{jk}-1}S(x,x^{m_j}t_{jk}(x)) = 0 \eqno (3.10a)$$ 
$$cx^{\alpha_j - m_jo_{jk}}|y - x^{m_j}t_{jk}(x)|\leq |\partial_y^{o_{jk}-1}S(x,y)| \leq Cx^{\alpha_j - m_jo_{jk}}|y - x^{m_j}t_{jk}(x)| \eqno (3.10b)$$
Since the terms of $S_{m_j}(x,y)$ are on the line $a + m_jb = \alpha_j$, the maximum possible order of a zero
of $S_{m_j}(1,y)$ is ${\alpha_j \over m_j}$. Hence $o_{jk} \leq {\alpha_j \over m_j}$ and the exponent in
$(3.10b)$ is a nonnegative number which we denote by $p$. 
Thus if we make the coordinate change $\phi_{jk}(x,y) = (x,y + x^{m_j}t_{jk}(x))$ and let $R_{jk}(x,y) = 
S\circ \phi_{jk}(x,y)$, then we have
$$\partial_y^{o_{jk}-1}R_{jk}(x,0) = 0 \eqno (3.11a)$$ 
$$cx^p|y| \leq |\partial_y^{o_{jk}-1}R_{jk}(x,y)| \leq Cx^p|y| \eqno (3.11b)$$
We now iterate the above algorithm to $R_{jk}(x,y)$ on $\phi_{jk}^{-1}(V_{jk})$. We first slice into two
pieces along the
$x$-axis. These two pieces are done the same way, so we focus our attention on the $y > 0$ piece which we
denote by $W$. We divide $W$ into pieces $T_{j'}'$, $U_{j'l'}'$, and $V_{j'k'}'$ exactly as done above. To
simplify notation, write $R(x,y) = R_{jk}(x,y)$, with the understanding that any subscripts
on $R$ really refer to subscripts on a fixed $R_{jk}$.

We will show that any positive zero $r_{j'k'}'$ of any $R_{m_j'}(1,y)$ has order at most $o_{jk} - 1$. 
Thus after at most $o_{jk}$ iterations there will be no more $V_{j'k'}'$ and we will be done. The fact 
that $R(x,y)$ is a smooth function of $x^{1 \over N}$ and $y$ rather than $x$ and $y$ does not cause 
any problems in future stages;
it just means after the next stage $N$ may be replaced by a large multiple of $N$. Also, there are 
no problems caused by the fact that the upper boundary of $\phi_{jk}^{-1}(V_{jk})$ is some curve 
$y = \xi x^{m_j} + $ higher order terms instead of $y = x^{\eta}$ as before; by shrinking $\xi$ if 
necessary we can ensure that this curve lies harmlessly inside one of the new $U_{j'l'}'$ whereupon the
only effect is to shrink this $U_{j'l'}'$ somewhat. 

So we turn our attention to showing that the order of any new positive zero $r_{j'k'}'$ of $R_{m_j'}(1,y)$ is 
at most $o_{jk} - 1$. For this, we will use $(3.11b)$. Note that such a zero occurs for $V_{j'k'}'$ of the 
form $\{(x,y): (r_{j'k'}' - \xi)x^{m_j'} < y < (r_{j'k'}' + \xi)x^{m_j'}\}$. The analogue of $(3.4)$ 
for $R(x,y)$ implies that on $V_{j'k'}'$ we have an expression 
$$\partial_{y}^{o_{jk}-1} R(x,y) = x^{\alpha_j'}\partial_{y}^{o_{jk}-1} R_{m_j'}(1,y) + 
x^{\alpha_j' + \zeta_j'} \partial_{y}^{o_{jk}-1}f(x,y) $$
$$+ O(|x|^M + |y|^{M - o_{jk}+1}) \eqno (3.12)$$
But by $(3.11b)$, along any curve $y = cx^{m_j'}$, the function $\partial_y^{o_{jk}-1}R(x,y)$ 
($= \partial_y^{o_{jk}-1}R_{jk}(x,y)$) will 
always have a zero of the same order $x^{p + m_j'}$ as $x \rightarrow 0$. Thus if the zero $r_{j'k'}'$
of $R_j'(1,y)$ were of order $o_{jk}$ or greater $(3.12)$ gives a contradiction: on the curve $y =
r_{j'k'}'x^{m_j'}$ the function $\partial_{y}^{o_{jk}-1} R(x,y)$ vanishes to order greater than
$\alpha_j'$, while on nearby curves $y = cx^{m_j'}$, $c \neq r_{j'k'}'$, it vanishes to order $\alpha_j'$. 
Thus we conclude that the order of the zero $r_{j'k'}'$ is at most $o_{jk} - 1$ and therefore that the
induction ends after finitely many steps. 

Lastly, part d) and the fact that each $g_i(x^{M'})$ is smooth for some positive integer $M'$ depending on $M$ 
follows by induction. Namely, if $V_{jk}$
comes from an edge $e$ with vertices $(a,b)$ and $(a',b')$, the coordinate change coming from that stage ofthe
induction is a smooth function of $x^{1 \over b' - b}$. In addition, each $o_{jk}$ is at most $b' - b$, and
the corresponding difference in $y$ coordinates will be at most $o_{jk}$ in future iterations. 
Also, the number of different $U_{jl}$ and $V_{jk}$ coming from that edge is bounded by twice the number of
possible zeroes of $S_e(1,y)$, or $2(b' - b)$. Since there are at most $o_{jk}$ iterations of the algorithm,
the result follows. This completes the proof of Theorem 3.1.

\noindent {\bf 4. The beginning of the proofs of Theorems 1.1 and 1.2; the first decomposition and 
preliminary lemmas.}  

We will prove Theorems 1.1 and 1.2 simultaneously; the proofs have much in common. Furthermore, in 
proving Theorem 1.1 we may assume that $S(x,y)$ has a degenerate critical point at $(0,0)$, as the 
nondegenerate (Morse) result can be read off from the explicit formulas (given for example in [G1]) 
for the leading term of the asymptotics for $M_{S,D}(\epsilon)$. Given this and that Theorem 1.2 excludes 
the nondegenerate case, the arguments of this section will always be under the assumption that $S(x,y)$ has
a degenerate critical point at $(0,0)$. Note also that we may assume that $S(x,y)$ is in superadapted
coordinates; a fixed coordinate change does not affect the statements of Theorem 1.1 or 1.2.

Note next that for a given vertex $v$ of $N(S)$ there is at most one value of $t$
for which the Taylor expansion of $S + tf$ at the origin does not have a $c_vx^v$ term. Thus there are
at most finitely many values of $t$ for which the Taylor expansion of $S + tf$ at the origin does not
contain a $c_vx^v$ term for each vertex $v$ of $N(S)$. In other words, other than for these values one
has $N(S) \subset N(S + tf)$. Thus excluding these values of $t$, in proving Theorems 1.1 and 1.2 we may 
assume that $N(S) \subset N(S + tf)$. 

Next, notice that in proving Theorem 1.1 it 
actually suffices to show $(1.4)$ holds for all sufficiently small $\epsilon
> 0$ and not for all $\epsilon$. The reason is as follows. Suppose for any sufficiently small 
disk $D$ centered at the origin one has $(1.3)$ for sufficiently small $\epsilon$. Fix one such 
neighborhood $D_0$. Then for all $\epsilon < \epsilon_0$ we have 
$$M_{S,D_0}(\epsilon) \leq C_S \epsilon^j |\ln \epsilon|^p$$
Since $M_{S,D}(\epsilon)$ is monotone in $D$, for all $D \subset D_0$, for $\epsilon < \epsilon_0$
one also has 
$$M_{S,D_0}(\epsilon) \leq C_S \epsilon^j |\ln \epsilon|^p$$
But by shrinking $D$ enough the inequality will also hold for ${1 \over 2} > \epsilon \geq \epsilon_0$. Thus if we
fix one such shrunken $D$, call it $D_1$, then $(1.4)$ will now hold for all ${1 \over 2} > \epsilon > 0$ for any
$D \subset D_1$. 

We now begin the proofs of Theorem 1.1 and 1.2. 
As indicated above, we can assume that $S(x,y)$ is in superadapted coordinates and that $N(S) 
\subset N(S + tf)$. Our goal will be to show that under the hypotheses of Theorem 1.1,  for all but 
finitely many $t$ $(1.4)$ holds for sufficiently small $\epsilon$, and that under the hypotheses of
Theorem 1.2, for all but finitely
many $t$ one has the estimate $M_{S + tf,D}(\epsilon) \leq C\epsilon^j|\ln \epsilon|^p$. Here $C$ may 
depend on $S, t, f$, and $D$.

We now fix a small disk $D$ centered at the origin. We divide $D$ into eight pieces through the $x$ and
$y$ axes and two curves $l_1$ and $l_2$ chosen as follows. If the bisectrix intersects the interior of 
a compact edge $e$ of $N(S)$ with equation $a + mb = \alpha$, then we choose $l_1$ and $l_2$ to be any 
two curves of the form $y = c|x|^m$ and $y = -c|x|^m$ so long as $c$ is not a zero of $S_e(1,y)$ or 
$S_e(-1,y)$. (Recall that the polynomial $S_e(x,y)$ is the sum of the terms $s_{ab}x^ay^b$ of $S(x,y)$'s 
Taylor expansion at the origin with $(a,b) \in e$). If
the bisectrix intersects $N(S)$ at a vertex $(d,d)$, we choose $l_1$ to be $y = |x|^m$ and $l_2$ to 
be $y = -|x|^m$, where $m$ is such that a line with equation $a + mb = \alpha$ intersects $N(S)$ at the
single point $(d,d)$. If the bisectrix intersects $N(S)$ in the interior
of the horizontal or vertical rays, switching the roles of the $x$ and $y$ axes if necessary we can 
assume it's the horizontal ray and $N(S)$'s lowest vertex is of the form $(c,d)$ for $c < d$. In this 
case we choose $l_1$ to be $y = |x|^m$ and $l_2$ to be $y = -|x|^m$ where a line with equation $a + mb = 
\alpha$ intersects $N(S)$ at the single point $(c,d)$. 

In the above fashion $D$ is divided into eight pieces $E_1,...E_8$. Clearly it suffices to show the 
desired bounds for each $M_{S + tf,E_i}(\epsilon)$. The argument for 
each $E_i$ is the same, so we will focus our attention on the piece from the upper right quadrant
between the $x$ axis and the curve $y = cx^m$ or $y = x^m$. Denote this piece by $E$. We now apply
the resolution of singularities algorithm of section 3 to $S(x,y)$, obtaining the corresponding sets
$D_i$. Define $D_i' = D_i \cap E$. Clearly it suffices to show the desired upper bounds for each
$M_{S + tf,D_i'}(\epsilon)$, which is what we will do.

Next, let $\phi_i$ be the maps of Theorem 3.1. Parts b) and c) of this theorem say that $S\circ 
\phi_i(x)$ is approximately
a monomial in the precise sense given there. Denote $\phi_i^{-1}D_i'$ by $F_i$, $S\circ \phi_i(x)$ by
$S_i(x)$, and $f \circ \phi_i(x)$ by $f_i(x)$. Since $\phi_i$ has determinant $\pm 1$, $M_{S + tf,D_i'}
(\epsilon) = M_{S_i + tf_i,F_i}(\epsilon)$. Thus our task is to prove good upper bounds for 
$M_{S_i + tf_i,F_i}(\epsilon)$. Also, in Theorem 1.1 the
smallness assumptions on the suprema of derivatives of $f$ are implied by corresponding smallness 
assumptions on the derivatives of $f_i$ (possibly with a different $\delta$), so in our future arguments
we may always assume the conditions on $f_i$ instead of $f$ without loss of generality. 

Analogous to above, excluding at most finitely many $t$ we can assume that $N(S_i) \subset N(R_i)$.
Also note that by Theorem 3.1, the set $F_i$ is a "curved triangle" in the sense that there are 
$h_i(x)$ and $H_i(x)$ such that for some $\eta_2 > \eta_1 > 0$ we have
$$\{(x,y): 0 < x < \eta_1,\,\,h_i(x) < y < H_i(x)\} \subset F_i \subset \{(x,y): 0 < x < \eta_2,\,\,
h_i(x) < y < H_i(x)\} $$
Here $h_i(x^{M'})$ and $H_i(x^{M'})$ are smooth for the $M$ given by the theorem. The function $h_i(x)$
may or may not be the zero function, and the Taylor expansion of $H_i(x)$ in fractional powers of $x$
has some nonvanishing initial term $A_ix^N_i$. Theorem 3.1 also says that we can let $(a_i,b_i)$ be
such that $S_i(x,y)$ is comparable to a multiple of $x^{a_i}y^{b_i}$ on $F_i$ in the sense of part b)
or c) of the theorem. The following lemma gives us various conditions satisfied by $(a_i,b_i)$, $N_i$,
and $N(S_i)$ what we will make use of.

\noindent {\bf Lemma 4.1.} 

\noindent {\bf a)} If the bisectrix intersects $N(S)$ in the interior of a compact edge, then if $S_i(x)$ 
is viewed as a function on the $x > 0$ part of a neighborhood of the origin, $S_i(x)$
is in superadapted coordinates with the same Newton distance $d$ that $N(S)$ has, and the bisectrix 
intersects $N(S_i)$ in the interior of a compact edge $e_i$. If the equation of this
edge is denoted by $a + m_ib = \alpha_i$, then $N_i \geq m_i$. Furthermore $b_i < d$ and $a_i > d$. 
The ordinate of the intersection of the bisectrix with the line of slope $-{1 \over N_i}$
containing $(a_i,b_i)$, given by ${a_i + N_ib_i \over  1 + N_i}$, is at most the Newton distance $d$ of
$N(S)$. 

\noindent {\bf b)} If the bisectrix intersects $N(S)$ at the vertex $(d,d)$, the same is true for
$N(S_i)$, and once again $S_i$ is in superadapted coordinates with Newton distance $d$. Either 
$(a_i,b_i) = (d,d)$, which happens for at least one $i$, or $b_i < d$ and $a_i > d$ like 
above. In the latter case we again have ${a_i + N_ib_i \over  1 + N_i} \leq d$. 

\noindent {\bf c)} If the bisectrix intersects $N(S)$ in the interior of one of the rays, then the same
is true for $N(S_i)$ and  again $S_i$ is in superadapted coordinates with Newton distance $d$. One of 
two things occurs. The first possibility is that $(a_i,b_i) = (c,d)$ for some $c < d$, the lower
boundary of $F_i$ is the $x$-axis, and part c) of Theorem 3.1 holds. The other possibility is that
$b_i < d$ and $a_i \geq d$. In either case, we have the strict inequality ${a_i + N_ib_i \over  1 + N_i} 
< d$.

\noindent {\bf Proof.} Recall that we divided a disk $D$ centered at the origin into
8 pieces, each of which after a coordinate change of the form $(x,y) \rightarrow (\pm x, \pm y)$ or 
$(\pm y, \pm x)$ becomes of the form $E = \{(x,y) \in D: 0 < y < cx^m\}$, and we are focusing our attention 
on $E$. In the new coordinates, $S(x,y)$ 
becomes a function which we denote by $S_0(x,y)$. Note that $S_0(x,y)$ is automatically still in 
superadapted coordinates. In the setting of part a) of this lemma, the bisectrix intersects $N(S_0)$ at the 
point $(d,d)$ which is in the interior of a compact edge $e$ with equation $a + mb = \alpha$ for some 
$\alpha$, $m$ as above. In the setting of part b), the intersection is still
$(d,d)$ which is now a vertex of $N(S_0)$, and in the setting of part c), the intersection is $(d,d)$ which
is in the interior of the horizontal ray. 

We now give some facts that are immediate consequences of the proof of the resolution of singularities 
algorithm of section 3, applied to $S_0(x,y)$ on $E$. First, each $\phi_i(x,y)$ is of 
the form $(x, \pm y - g_i(x))$, where $g_i(x^K)$ is a smooth function for some $K$. Next, if $g_i(x)$ is
not the zero function then the Taylor 
expansion of $g_i(x)$ in fractional powers of $x$ has initial term $p_ix^{l_i}$ where $l_i \geq m$ is such 
that $N(S_0)$ has an edge with equation $a + l_i b = \alpha_i$ for some $i$. The definition of $E$ is such
that $(d,d)$ is either the upper vertex of this edge, or the edge lies entirely below $(d,d)$. The number 
$N_i$, defined such that the upper boundary of $F_i$ has equation $y = q_i x^{N_i} + $ higher order terms, 
satisfies $N_i \geq m$ and the algorithm ensures that $(a_i,b_i)$ is a vertex of $N(S_i)$. The definition of 
$E$ is such that in the settings of part a) or b) of this lemma, either $(a_i,b_i) = (d,d)$ and $N_i = m$,
or $(a_i,b_i)$ is the lower vertex of a compact edge of $N(S_i)$ of slope $-{1 \over N_i}$.
In the latter case either $(d,d)$ is the upper vertex of
the edge or it lies entirely below $(d,d)$.

We now proceed with the proof, starting with part a). If $l_i > m$, then the coordinate change $\phi_i(x,y)$
will not affect $(S_0)_e(x,y)$, the sum of the terms of $S_0(x,y)$'s Taylor series at $(0,0)$ corresponding to 
the edge $a + mb = \alpha$. Hence the resulting function $S_i(x,y)$ 
will be in superadapted coordinates like before. On the other hand, if $l_i = m$, then $(S_0)_e(x,y)$ becomes
$(S_0)_e(x, \pm y \ - p_ix^m)$ after the coordinate change. Hence $(S_0)_e(1,y)$ becomes $(S_0)_e(1, \pm y - p_i)$ and 
$(S_0)_e(-1,y)$
becomes $(S_0)_e(-1, \pm y - p_i)$. Shifting the $y$ variable does not change the fact that the definition of 
superadapted coordinates holds; the condition is that these polynomials have zeroes of order
less than $d$. Hence in the $l_i = m$ case we are in superadapted coordinates as well. In either case,
the bisectrix still intersects $N(S_i)$ at $(d,d)$, which is in the interior of a compact edge $e_i$ with 
equation $a + mb = \alpha$ for some $\alpha$. As a result $(a_i,b_i)$, being a vertex of $N(S_i)$ lying
below $(d,d)$, satisfies $a_i > d$ and $b_i < d$ as needed. Using the last paragraph, $(a_i,b_i)$ is the 
lower vertex of a compact edge $e_i$ of $N(S_i)$ with slope $-{1 \over N_i}$ which either contains $(d,d)$ 
or is below the edge containing $(d,d)$.
Thus the intersection of the bisectrix with the line containing $e_i$ is at $(d,d)$
or below. But $(a_i,b_i)$ is on this edge, so the intersection occurs at $({a_i + N_ib_i \over 1 + N_i},
{a_i + N_ib_i \over 1 + N_i})$. Hence ${a_i + N_ib_i \over 1 + N_i} \leq d$ as needed. This completes the
proof of part a).

Next we suppose we are in the setting of part b). In this case, the initial term of $g_i(x)$ is 
$p_ix^{l_i}$ for some $l_i > m$ since $N(S_0)$ has no edge with equation of the form $a + mb = \alpha$.
As a result, the coordinate change $\phi_i(x,y)$ will not alter the fact that the bisectrix intersects
the Newton polygon at $(d,d)$. Furthermore, one will still be in superadapted coordinates; if $e$ is an 
edge of $N(S_0)$ containing $(d,d)$, then either $(S_0)_e(x,y)$ is unchanged by the coordinate change, or
$(S_0)_e(1,y)$ becomes $(S_i)_e(1,y) = (S_0)_e(1,\pm y - p_i)$ and $(S_0)_e(-1,y)$ becomes $S_i(-1,y) = 
(S_0)_e(-1,\pm y - p_i)$ like above. In either event $S_i(x,y)$ will still be in superadapted coordinates. 
By the last paragraph, $(a_i,b_i)$ is either $(d,d)$ or a lower vertex. In the former case, ${a_i + N_ib_i 
\over 1 + N_i} =  d$, and
in the latter case, exactly as in part a) we have ${a_i + N_ib_i \over 1 + N_i} \leq d$ as needed.
This completes the proof of part b).

Suppose we are in the setting of part c). Then the bisectrix intersects $N(S_0)$ either in the interior of
the horizontal or vertical ray. Suppose it is the horizontal ray. Then $N(S_0)$ has a vertex $(c,d)$ for
some $c < d$. In this case no further subdivisions are necessary; we can take $S_i = S_0$ and let $F_i$ be
all of $\{(x,y) \in D: 0 < y < x^m\}$. Part c) of Theorem 3.1 automatically holds. So suppose the bisectrix
intersects $N(S_0)$ in the interior of the vertical ray. In this case, the highest vertex of $N(S_0)$ is
$(d,c)$ for some $c < d$. Since any coordinate change $\phi_i$ fixes the highest vertex of $N(S_0)$, the
highest vertex of $N(S_i)$ is $(d,c)$ as well. Thus $S_i$ is in superadapted coordinates with the bisectrix
intersecting $N(S_i)$ inside its vertical ray. Since $(a_i,b_i)$ is either the vertex $(d,c)$ or a lower 
one, we have $a_i \geq d$ and $b_i \leq c < d$. Lastly, since the bisectrix intersects $N(S_i)$ inside the
vertical ray, the ordinate of the intersection of any 
supporting line of $N(S_i)$ containing $(a_i,b_i)$ with $N(S_i)$ is less than $d$. Hence ${a_i + N_ib_i
\over 1 + N_i} < d$ and we are done. 

\noindent {\bf 5. The Main Argument.} 

We work in the setting of section 4. For some fixed value of $t$, let $R_i = S_i + tf_i$. We will proceed 
along the lines outlined in section 1,
dividing a given $F_i$ into finitely many pieces and show that, excluding finitely many values of
$t$, the contribution to $M_{R_i,F_i}(\epsilon)$
of each piece satisfies the upper bounds given by Theorem 1.1 or 1.2. We start this as
follows. Since $S_i(x^{M'},y)$ and $f_i(x^{M'},y)$ are 
smooth functions, where $M$ is as in Theorem 3.1, by Corollary 3.2, we can apply the resolution of
singularities algorithm to $R_i(x,y)$. We now do so, but focus our attention only on the first 
stage of the algorithm, dividing the upper right quadrant into the sets called $T_j$ and $U_j$ in
the proof. To highlight their dependence on $i$ here, we refer to them as $T_{ij}$ and $U_{ij}$ here.
Clearly, we need only consider those $T_{ij}$ and $U_{ij}$ that intersect $F_i$. Each $U_{ij}$ 
corresponds to an edge of $N(R_i)$ whose equation we write as $a + m_{ij}b = \alpha_{ij}$, while each
$T_{ij}$ corresponds to a vertex of $N(R_i)$. 

If $N(R_i)$ has a compact edge whose upper vertex is on or below the line $y = b_i$ and has slope greater 
than $-{1 \over N_i}$, let $e_{ij'}$ be the uppermost amongst all such edges. If
$e_{ij'}$ exists, let $G_i$ denote the union of $U_{ij'}$ with 
any $T_{ij}$ and $U_{ij}$ corresponding to edges and vertices of $N(R_i)$ below $e_{ij'}$. If $e_{ij'}$ 
does not exist, simply define $G_i$ to be the lowest $T_{ij}$. We will now find upper bounds for
$M_{R_i,G_i}(\epsilon)$ that are as good as needed for Theorem 1.1 or 1.2.

\noindent {\bf Lemma 5.1.} Under the assumptions of Theorem 1.1, we have $M_{R_i,G_i}(\epsilon) \leq
C_S \epsilon^j |\ln \epsilon|^p$, and under the assumptions of Theorem 1.2 we have $M_{R_i,G_i}
(\epsilon) \leq A \epsilon^j |\ln \epsilon|^p$ for some constant $A$.

\noindent {\bf Proof.} Let $(a',b')$ be the uppermost vertex of the union of all edges and vertices of
$N(R_i)$ that correspond to a $T_{ij}$ or $U_{ij}$ included in $G_i$. If $G_i$ consists solely of the
lowest $T_{ij}$, let $(a',b')$ be the lowest vertex of $N(R_i)$. Note that by definition of $G_i$,
we have $b' \leq b_i$.
Write the Taylor expansion of $R_i(x,y)$ at the origin as $\sum_{cd}r_{cd}x^cy^d$, so that $r_{a'b'}
x^{a'}y^{b'}$ denotes the term corresponding to the vertex $(a',b')$. Note that this Taylor expansion
may contain fractional powers of $x$, but not $y$. By Lemma 2.4 c), on $G_i$ we have 
$$|\partial_y^{b'} R_i(x,y)| > {b'! \over 2}|r_{a'b'} x^{a'}| \eqno (5.1)$$
Thus by Lemma 2.2 we have
$$M_{R_i,G_i}(\epsilon) = |\{(x,y) \in G_i: |R_i(x,y)| < \epsilon\}| $$
$$\leq 4 |\{(x,y) \in G_i: {1 \over 2}|r_{a'b'}x^{a'}y^{b'}| < \epsilon\}| \eqno (5.2)$$
The right-hand side of $(5.2)$ may be estimated using Lemma 2.3. We break into cases, depending on whether
$b' > a'$, $b' = a'$, or $b' < a'$. If $b' > a'$, the lemma says that
the right-hand side of $(5.2)$ is bounded by $Cr^{\eta}\epsilon^{1 \over b'}$ where $r$ is the radius
of $D$. Since $b' \leq b_i \leq d$, $d$ the Newton distance of $S$, this can be made less than
$\epsilon^{1 \over d}$ by shrinking the radius of original disk $D$ enough. Since the growth index of 
$S(x,y)$ is given by ${1 \over d}$ in superadapted coordinates, this is at least as good as the 
estimate we need.

\noindent If $b' = a'$, then Lemma 2.3 tells us that 
$$|\{(x,y) \in G_i: {1 \over 2}|r_{b'b'}x^{b'}y^{b'}| < \epsilon\}| < C\epsilon^{1 \over b'}
|\ln \epsilon| \eqno (5.3)$$
Here the constant $C$ depends on lower bounds for $|r_{b'b'}|$ as well as the set $A_1$.
Since $b' \leq b_i \leq d$, this is better than the estimate we seek unless $b' = b_i = d$. So we suppose
$b' = b_i = d$. Since $N(S_i) \subset N(R_i)$, we must have $a_i \geq a' = d$. By Theorem 4.1, the only
way one can have $a_i \geq d$ and $b_i = d$ is for $(a_i,b_i) = (d,d)$. Hence we have $(a_i,b_i) = 
(a',b') = (d,d)$. In view of Theorem 3.1c), $F_i$ contains a set of the 
form $\{(x,y): 0 < x < \eta,\,\,x^{m_1} < y < x^{m_2}\}$ on which $S_i(x,y) \sim x^dy^d$, and so by Lemma 2.3b)
one has that $M_{S_i,G_i}(\epsilon) > c\epsilon^{1 \over d}|\ln \epsilon|$. Hence we must be in the situation
where the growth index ${1 \over d}$ of $S(x,y)$ has multiplicity 1, and $(5.2)$ and $(5.3)$ give 
the desired estimate for Theorem 1.2. Next, note that
$$|\{(x,y) \in G_i: {1 \over 2}|r_{dd}x^{d}y^{d}| < \epsilon\}| = |\{(x,y) \in G_i: x^{d}y^{d} <
{2\epsilon \over r_{dd}} \}| \eqno (5.4)$$
In the setting of Theorem 1.1, by making $\delta$ sufficiently small we can
ensure that $r_{dd} <  2s_{dd}$, where $s_{dd}x^dy^d$ denotes the term of the
Taylor expansion of $S_i(x,y)$. Hence, using Lemma 2.3b) on $(5.2)$ and $(5.4)$, we see that $M_{R_i,G_i}
(\epsilon) < C_S\epsilon^{1 \over d} |\ln \epsilon|$, the desired estimate for Theorem 1.1.

We now turn to the case where $b' < a'$. The definition of $G_i$ is such that $G_i$ is 
contained in a set of the form $\{(x,y): 0 < x < \eta,\,\,0 < y < x^m\}$ for some $m > N_i$.
We apply part c) of Lemma 2.3, which says that 
$$|\{(x,y) \in G_i: {1 \over 2}|r_{a'b'}x^{a'}y^{b'}| < \epsilon\}| \leq C\epsilon^{m + 1 \over a'
+ mb'} \eqno (5.5a)$$
In view of $(5.2)$, we have
$$M_{R_i,G_i}(\epsilon) = |\{(x,y) \in G_i: |R_i(x,y)| < \epsilon\}| < C'\epsilon^{m + 1 \over a' + mb'}
\eqno (5.5b)$$
Like above, the constant $C'$ depends on lower bounds for $R_{a'b'}$ (as well as $(a',b')$).
Note that ${a' + mb' \over m + 1}$ is the ordinate of the intersection of the line of slope
$-{1 \over m}$ containing $(a',b')$ with the bisectrix. Since $a' > b'$ and $m > N_i$, this is strictly
less than the corresponding ordinate for the line of slope $-{1 \over N_i}$ containing $(a',b')$.
Since $N(R_i)$ contains $N(S_i)$, this will be at most the ordinate of the intersection of the bisectrix
with line of the same slope $-{1 \over N_i}$ containing $(a_i,b_i)$, given by ${a_i + b_iN_i \over 1 + 
N_i}$. This is at most $d$ by Lemma 4.1. Hence by $(5.5b)$, we have $M_{R_i,G_i}(\epsilon) < C\epsilon^{{1 
\over d} + \eta}$ for some positive $\eta$, a better estimate than we need and we are done. 

\noindent The next step is to prove upper bounds of Theorem 1.1 and 1.2 for the remaining $M_{R_i,T_{ij}}
(\epsilon)$:

\noindent {\bf The upper bounds for $M_{R_i,T_{ij}}(\epsilon)$ when $T_{ij}$ is not a subset of $G_i$.} 

Each such $T_{ij}$ corresponds to some vertex of $N(R_i)$, which we denote by $(p,q)$.  
$T_{ij}$ is typically of the form $\{(x,y) \in F_i: {1 \over \xi}x^{1 \over m_2} < y < \xi x^{m_1}\}$. 
It's possible that an uppermost $T_{ij}$ is some proper subset of such a set, but the following argument
works for that situation too. Also, although
we only need to bound $M_{R_i,T_{ij}}(\epsilon)$ for $T_{ij}$ intersecting $F_i$, the following argument
works for all $T_{ij}$. 
We will analyze $R_i(x,y)$ on $T_{ij}$ similarly to the way $R(x,y)$ was analyzed in $(2.8)$. As above,
we write the Taylor series of $R_i(x,y)$ at the origin as $\sum_{ab} r_{ab}x^ay^b$. We similarly write 
the Taylor series of $S_i(x,y)$ at the origin as $\sum_{ab} s_{ab}x^ay^b$. Note that 
$$S_i(x,y) -  s_{pq}x^py^q= \sum_{p \leq a < M,\,\,q \leq b < M,\,\,(a,b) \neq (p,q)}s_{ab}
x^ay^b$$
$$ + \sum_{a < p,\,\,q < b < M ,\,\,a + m_1b \geq \alpha_1}s_{ab}x^ay^b + 
\sum_{p < a < M,\,\,b < q,\,\,am_2 + b \geq \alpha_2}s_{ab}x^ay^b + E_M(x,y) \eqno (5.6)$$
Here $a + m_1b = \alpha_1$ and $m_2a + b = \alpha_2$ are equations of the edges $e_1$ and $e_2$ of
$N(R_i)$ meeting at $(p,q)$. In the case that $e_1$ doesn't exist, the second sum above disappears, 
and if $e_2$ doesn't exist then the third sum disappears.

First suppose $(p,q)$ is not a vertex of $N(S)$, so that $s_{pq} = 0$. By definition, $r_{pq}$ is 
nonzero. As in $(2.8)$, for any $\delta > 0$ we can make each of the terms of $(5.6)$ bounded by $\delta 
|r_{pq}|x^py^q$ by choosing $\xi$ small enough. In particular, we can make the
absolute value of the entire right hand side of $(5.6)$ less than ${1 \over 4}|r_{pq}|x^py^q$ on  
$T_{ij}$. Since $s_{pq} = 0$ this means $|S_i(x,y)| < {1 \over 4}|r_{pq}|x^py^q$ on $T_{ij}$. Similarly, 
by choosing $\xi$ small enough, we can assume that the
right-hand side of the analogue to $(5.6)$ with $S_i$ replaced by $R_i$ is also less than 
${1 \over 4}|r_{pq}|x^py^q$ on $T_{ij}$. Hence on $T_{ij}$ we have
$$|R_i(x,y)| > {3 \over 4} |r_{pq}|x^py^q > 3|S_i(x,y)| \eqno (5.7)$$
As a result, we have
$$|\{(x,y) \in T_{ij}: |R_i(x,y)| < \epsilon\}| \leq |\{(x,y) \in T_{ij}: |S(x,y)| < 3\epsilon\}|$$ 
$$ \leq C_S \epsilon^j|\ln \epsilon|^p$$
This is at least as strong as the estimate we need.

Now suppose that $s_{pq} \neq 0$. Since $(p,q)$ is a vertex of $N(R_i)$ and $N(S) \subset N(R_i)$, this
means $(p,q)$ is a vertex of $N(S)$ as well. Like above, by shrinking $\xi$ enough one can assume that
the right-hand side of $(5.6)$ and its analogue with $S_i$ replaced by $R_i$ are less than 
${1 \over 2}|s_{pq}|$ and ${1 \over 2}|r_{pq}|$ respectively. As a result, on $T_{ij}$ we have
$$|S_i(x,y)| < {3 \over 2}|s_{pq}|x^py^q,\,\,\,|R_i(x,y)| > {1 \over 2}|r_{pq}|x^py^q$$
Thus we have
$$|\{(x,y) \in T_{ij}: |R_i(x,y)| < \epsilon\}| < |\{(x,y) \in T_{ij}: |S_i(x,y)| < 3
{|s_{pq}| \over |r_{pq}|}\epsilon\}| \eqno (5.8)$$
The right-hand side of $(5.8)$ is at most $C \epsilon^j|\ln \epsilon|^p$,
the desired estimate in the setting of Theorem 1.2. As for Theorem 1.1, we may assume the 
$\delta$ of that Theorem is small enough so that ${|s_{pq}| \over |r_{pq}|} < 2$. In this case, we have
$$|\{(x,y) \in T_{ij}: |S_i(x,y)| < 3 {|r_{pq}| \over |s_{pq}|}\epsilon\}| < 
|\{(x,y) \in T_{ij}: |S_i(x,y)| < 6 \epsilon\}|$$
$$< C_S \epsilon^j|\ln \epsilon|^p \eqno (5.9)$$
Combining $(5.8)$ and $(5.9)$ gives the desired estimates for Theorem 1.1 and we are done.

It remains to bound $M_{R_i,{U_{ij}}}(\epsilon)$ for the $U_{ij}$ intersecting $F_i$ that are not contained 
in $G_i$. Each $U_{ij}$ corresponds
to an edge $e_{ij}$ of $N(R_i)$. There are only finitely many possible such $e_{ij}$ for 
{\it any} $R_i$; that is, there are only finitely many pairs of vertices that can be the endpoints
of an edge corresponding to some such $U_{ij}$ for any possible $S_i + tf_i$, regardless of what $f_i$
is. To see why this is true, write the equation of $e_{ij}$ as $a + m_{ij}b = \alpha_{ij}$, and its upper
vertex as $(a,b)$. Since $U_{ij}$ us a subset of $F_i$, we automatically have $m_{ij} \leq N_i$. We 
separately consider the cases $b < b_i, b= b_i$, and $b > b_i$, and show in each case that there are
finitely many possibilities for $e_{ij}$. 

We start with the case where $b < b_i$. In this case $G_i$ takes care of all $U_{ij}$ with $m_{ij} 
> N_i$, so we are left with the case when $m_{ij} = N_i$. There are only finitely many possibilities
with $a \leq a_i$, so we may assume that $a > a_i$. Note that since $(a_i,b_i)$ is in $N(S_i)$, it is
also in $N(R_i)$. Since $a > a_i$ this means $(a,b)$ can not be on a the vertical ray of $N(R_i)$.
Instead it is the lower vertex of a compact edge $e'$ of $N(R_i)$ of slope at least $ {b - b_i \over 
a - a_i} \geq -{1 \over a - a_i}$. Since the slope of $e_{ij}$ is $-{1 \over N_i}$ and $e_{ij}$ lies 
below $e'$, we must have that $-{1 \over N_i} < -{1 \over a - a_i}$. In other words we have $a < a_i 
+ N_i$. Since it also true that $b < b_i$, there are finitely many possibilities for this to occur.

Next, consider the case where $b = b_i$. Here since $N(R_i)$ contains $N(S_i)$ and $(a_i,b_i)$ is a 
vertex of $N(S_i)$, $a \leq a_i$ and there are finitely many possibilities for $a$. Since once again
$G_i$ takes care of all $U_{ij}$ with $m_{ij} > N_i$, we have $m_{ij} = N_i$ and once again there are 
a finite number of possibilities for $e_{ij}$.

Lastly, we consider the situation where $b > b_i$. Then since $(a,b)$ lies above the vertex $(a_i,b_i)$
of $N(S_i)$, we have $a < a_i$. Since $N(R_i)$ contains $N(S_i)$, $(a_i,b_i) \in N(R_i)$. Hence the
slope of of $e_{ij}$ is at most ${b_i - b \over a_i - a} \leq b_i - b$. It is also greater than or equal
to $-{1 \over N_i}$ because $U_{ij} \subset F_i$. Thus we have $-{1 \over N_i} \leq b_i - b$ or $b \leq b_i
+ {1 \over N_i}$. 
Coupled with the condition that $a < a_i$, once again we have finitely many possibilities for $(a,b)$
and we are done. 

Thus for our future arguments, it suffices to fix a single $e_{ij}$ and prove the desired upper bounds
for the $U_{ij}$ associated with $e_{ij}$. Recall that for a fixed $x$ the vertical cross-section of
$U_{ij}$ has width $({1 \over \xi} - \xi)x^{m_{ij}}$. If one now applies the next step of the
resolution of singularities algorithm of Theorem 3.1 (to $R_i(x,y)$), one divides $U_{ij}$ into 
pieces $V_{ijk}$ and $W_{ijk}$, where each $V_{ijk}$ is of the form $\{(x,y): (r_{ijk} - \xi)x^{m_{ij}} 
< y < (r_{ijk} + \xi)x^{m_{ij}}$ for a root $r_{ijk}$ of the polynomial $(S_i)_{e_{ij}}(1,y)$. On each
$W_{ijk}$ one has $R_i(x,y) \sim x^{\alpha_{ij}}$ where $e_{ij}$ has equation $a + m_{ij}b = \alpha_{ij}$.
We also need to split each $V_{ijk}$ into two pieces $V_{ijk}^1$, and $V_{ijk}^2$, where 
$V_{ijk}^2$ is the portion where $|y - r_{ijk}x^{m_{ij}}| < x^{m_{ij} + \eta}$, and $V_{ijk}^1$ is the
rest. $\eta$ here is a small positive constant which must be sufficiently small for our arguments to 
work. Note that since there are at most boundedly many $r_{ijk}$ for a given $e_{ij}$ and boundedly 
many $e_{ij}$ for our fixed $S(x,y)$, the total number of $V_{ijk}^1, V_{ijk}^2$, and $U_{ijk}$ is 
uniformly bounded given $S(x,y)$. 

Our analysis will proceed as follows. We will first show that each $M_{R_i, W_{ijk}}(\epsilon)$ 
satisfies the bounds required of Theorems 1.1 or 1.2. Then we will show that if $\eta > 0$ is 
sufficiently small, a small additional argument will show that each $M_{R_i, V_{ijk}^1}(\epsilon)$
also satisfies the bounds. Afterwards, a separate argument will be used to show that for each $V_{ijk}^2$ 
intersecting
$F_i$, $M_{R_i, V_{ijk}^2}(\epsilon)$ satisfies bounds better than the ones we need. We will do this by 
applying the full resolution of singularities algorithm on $f_i(x,y)$ on each $V_{ijk}^2$.
We will obtain the corresponding sets called $D_i$ in Theorem 3.1, and show that for each $D$ amongst them
$M_{R_i, D}(\epsilon)$ satisfies better estimates than what we need. 

\noindent {\bf The analysis of $M_{R_i, W_{ijk}}(\epsilon)$.}

\noindent Since $N(R_i)$ contains $N(S_i)$ and $e_{ij}$ has equation $a + m_{ij}b = \alpha_{ij}$, 
$(a_i,b_i)$ is on or above the line containing $e_{ij}$ and hence $a_i + m_{ij}b_i \geq 
\alpha_{ij}$. Since on $W_{ijk}$ we have $R_i(x,y) \sim x^{\alpha_{ij}}$, we have
$R_i(x,y) \geq C_1x^{a_i + m_{ij}b_i}$. Since for fixed $x$ the set $W_{ijk}$ has cross-section
$\sim x^{m_{ij}}$, we conclude that 
$$|\{(x,y) \in W_{ijk}: |R_i(x,y)| < \epsilon| < C_2\epsilon^{1 + m_{ij} \over a_i + m_{ij}b_i}
\eqno (5.10a)$$
If $(a_i,b_i)$ is strictly above the line, then for some $\zeta > 0$ we have the even better estimate
$$|\{(x,y) \in W_{ijk}: |R_i(x,y)| < \epsilon| < C_3\epsilon^{{1 + m_{ij} \over a_i + m_{ij}b_i} + \zeta}
\eqno (5.10b)$$
For now at least, we have no information concerning the constants $C_2,C_3$ of $(5.10a)-(5.10b)$. Note that
the exponent ${1 + m_{ij} \over a_i + m_{ij}b_i}$ is the reciprocal of the ordinate of the intersection
of the bisectrix with the line of slope $-{1 \over m_{ij}}$ containing $(a_i,b_i)$. Also note that since
$S_i(x,y) \sim x^{a_i + m_{ij}b_i}$ on $W_{ijk}$, we have 
$$|\{(x,y) \in W_{ijk}: |S_i(x,y)| < \epsilon| > C_4\epsilon^{1 + m_{ij} \over a_i + m_{ij}b_i}
\eqno (5.11)$$
Hence one has
$$|\{(x,y) \in W_{ijk}: |R_i(x,y)| < \epsilon| < C_5 |\{(x,y) \in W_{ijk}: |S_i(x,y)| < \epsilon|$$
$$ < C_5 M_{S,D}(\epsilon)$$
$$\leq C_6\epsilon^{1 \over d}|\ln \epsilon|^p \eqno (5.12)$$
This gives the desired estimates for Theorem 1.2, which are also the desired estimates for 
Theorem 1.1 other than the constants, which will take some more work and which we now focus our attention on.

Note that if $(a_i,b_i)$ lies strictly above the line containing $e_{ij}$, the added $\zeta$ in $(5.10b)$ 
gives us any constant $C_6$ we want for $\epsilon$ small enough, thereby implying the estimate needed for 
Theorem 1.1. Thus we need only consider
the case where $(a_i,b_i)$ is actually on the line containing $e_{ij}$. Furthermore, in any situation
in which ${1 + m_{ij} \over a_i + m_{ij}b_i}$ is strictly greater than ${1 \over d}$, we could once 
again make $C_6$ arbitrarily small. So it makes sense to analyze when we have
do not have this strict inequality; we will see momentarily that this only happens in special situations.

Consider the case when the bisectrix intersects the interior of a compact edge of $N(S_i)$. In
this case, by Lemma 4.1a) we have $b_i < a_i$. Hence if $m_{ij}$ is greater than the minimal possible
value on $F_i$, given by $N_i$, then we have ${1 + m_{ij} \over a_i + m_{ij}b_i} > {1 + N_i \over a_i 
+ N_ib_i}$, which in turn is at least ${1 \over d}$. Hence when the bisectrix intersects the interior
of a compact edge of $N(S_i)$, equality can only occur if $m_{ij} = N_i$.

Next, consider the case when the bisectrix intersects the vertex $(d,d)$ of $N(S_i)$. In this case,
by Lemma 4.1b) either $(a_i,b_i) = (d,d)$ or $a_i > d, b_i < d$, and ${1 + N_i \over a_i 
+ N_ib_i} \geq {1 \over d}$. In the latter case, if $m_{ij} > N_i$ then like above we have ${1 + m_{ij}
\over a_i + m_{ij}b_i} > {1 + N_i \over a_i + N_ib_i} \geq {1 \over d}$, and equality can occur only 
when $m_{ij} = N_i$.

Lastly, consider the case where the bisectrix intersects $N(S_i)$ in the interior of one of the rays.
In this case Lemma 4.1c) applies and either $(a_i,b_i) = (c,d)$ for some $c < d$, or $(a_i,b_i)$
satisfies $a_i \geq d$, $b_i < d$ and ${1 + N_i \over a_i + N_ib_i} > {1 \over d}$. In the latter case
since $m_{ij} \geq N_i$ we automatically have ${1 + m_{ij} \over a_i + m_{ij}b_i} > {1 \over d}$ and 
equality does not occur. In the former case, since $c < d$ we also have ${1 + m_{ij} \over c + m_{ij}d} > 
{1 \over d}$ and equality also does not occur.

In summary, in bounding $M_{R_i, W_{ijk}}(\epsilon)$ we have already covered all possible cases except 
the following two situations. First, we can have $m_{ij} = N_i$, $a_i > d$, and $b_i < d$. Secondly, the 
bisectrix may intersect $N(S_i)$ at the vertex $(d,d)$ with $(a_i,b_i) = (d,d)$. Furthermore, as mentioned 
above we only need to consider the situation when 
$(a_i,b_i)$ is on the line containing the edge $e_{ij}$. The argument we will use for these two
situations actually will give the needed bounds for all of $M_{R_i,V_{ijk}}(\epsilon)$ in those situations
as well. In fact, recalling that $U_{ij} = \cup_k V_{ijk} \cup \cup_k W_{ijk}$, what we will do is prove 
upper bounds for all of $M_{R_i,U_{ij}}(\epsilon)$ in both situations.

\noindent {\bf Bounding $M_{R_i, U_{ij}}(\epsilon)$ in the two exceptional situations.} 

We start with the first situation, where $m_{ij} = N_i$, $a_i > d$, $b_i < d$, and $(a_i,b_i)$ is on the
line containing $e_{ij}$. We write the Taylor expansions of $R_i(x,y)$ and $S_i(x,y)$ at the origin
as $\sum_{a,b} r_{ab}x^ay^b$ and $\sum_{a,b} s_{ab}x^ay^b$. Here the $b$'s are all integers but $a$
may be a nonintegral positive rational number. We write
$$R_i(x,y) = \sum_{a + N_ib = \alpha_{ij}} r_{ab}x^ay^b + \sum_{a + N_ib > \alpha_{ij},\,\,a,b < M}
r_{ab}x^ay^b + O(|x|^M + |y|^M) \eqno (5.13a)$$
$$S_i(x,y) = \sum_{a + N_ib = \alpha_{ij}} s_{ab}x^ay^b + \sum_{a + N_ib > \alpha_{ij},\,\,a,b < M}
s_{ab}x^ay^b + O(|x|^M + |y|^M) \eqno (5.13b)$$
Here $M$ is some large integer and $\alpha_{ij}$ is such that $e_{ij}$ has equation $a + N_ib = 
\alpha_{ij}$. Since for fixed $x$, the $y$ cross-section of $U_{ij}$ is contained in some 
$[c_{ij}x^{N_i},A_ix^{N_i}]$ where $c_{ij}, A_i > 0$ and $A_i$ depends only on $S_i$, it makes sense to 
look at $R_i(x,x^{N_i}y)$ and $S_i(x,x^{N_i}y)$, given by
$$R_i(x,x^{N_i}y) = x^{\alpha_{ij}}(R_i)_{e_{ij}}(1,y) + x^{\alpha_{ij}+\zeta}P(x,y) + O(|x|^{M'})
\eqno (5.14a)$$
$$S_i(x,x^{N_i}y) = x^{\alpha_{ij}}(S_i)_{e_{ij}}(1,y) + x^{\alpha_{ij}+\zeta}Q(x,y) + O(|x|^{M'})
\eqno (5.14b)$$
Here $\zeta > 0$, $P$ and $Q$ are polynomials in $y$ and a fractional power of $x$, and $M'$ is a large
integer that grows linearly with $M$. Similarly, we may look at the $b_i$th $y$ derivatives of $R_i$ and 
$S_i$ and get the following expressions
$$(\partial_y^{b_i} R_i)(x,x^{N_i}y) = x^{\alpha_{ij} - b_iN_i}\partial_y^{b_i}((R_i)_{e_{ij}}(1,y)) + 
x^{\alpha_{ij} - b_iN_i + \zeta}\tilde{P}(x,y) + O(|x|^{M'}) \eqno (5.15a)$$
$$(\partial_y^{b_i} S_i)(x,x^{N_i}y) = x^{\alpha_{ij} - b_iN_i}\partial_y^{b_i}((S_i)_{e_{ij}}(1,y))+ 
x^{\alpha_{ij} - b_iN_i + \zeta}\tilde{Q}(x,y) + O(|x|^{M'}) \eqno (5.15b)$$
Part b) or c) of Theorem 3.1 says that $\partial_y^{b_i} S_i(x,y) > Cx^{a_i} =  Cx^{\alpha_{ij} - 
b_iN_i}$ on the $F_i$ which contains the set $U_{ij}$ under consideration. Hence letting $x \rightarrow 0$ in
$(5.15b)$, we may conclude that $\partial_y^{b_i}((S_i)_{e_{ij}}(1,y))$ has no zeroes on $[0,A_i]$ 
Hence there is a $C_S$ such that for the $(x,y)$ being considered we have
$$|\partial_y^{b_i} ((S_i)_{e_{ij}}(1,y))x^{\alpha_{ij} - b_iN_i}| > C_Sx^{\alpha_{ij} - b_iN_i}
\eqno (5.16)$$
Suppose now we are working under the hypotheses of Theorem 1.1. Then by shrinking $\delta$ sufficiently
if necessary, we can assume each coefficient of each term of $(\partial_y^{b_i} f_i)_{e_{ij}}(1,y)$ is as 
small as we want. In particular, we can assume that on $[c_{ij},A_i]$ we have
$$|\partial_y^{b_i} ((f_i)_{e_{ij}}(1,y))x^{\alpha_{ij} - b_iN_i}| < {1 \over 2} C_Sx^{\alpha_{ij}
- b_iN_i} \eqno (5.17a)$$
Given that $R_i = S_i + tf_i$, for $|t| < 1$ this means that 
$$|\partial_y^{b_i}(( R_i)_{e_{ij}}(1,y))x^{\alpha_{ij} - b_iN_i}| > {1 \over 2} C_Sx^{\alpha_{ij} - b_iN_i}
\eqno (5.17b)$$
In view of $(5.15a)$, this means that if $x$ is sufficiently small, which we may assume, then we have
$$|(\partial_y^{b_i} R_i)(x,x^{N_i}y)| >  {1 \over 4}C_Sx^{\alpha_{ij} - b_iN_i} \eqno (5.18a)$$
Translating this back into the original coordinates, this means that on $U_{ij}$ we have
$$|\partial_y^{b_i} R_i(x,y)| >  {1 \over 4}C_Sx^{\alpha_{ij} - b_iN_i} \eqno (5.18b)$$
Hence by Lemma 2.2, for a given $x$ we have
$$|\{y: |R_i(x,y)| < \epsilon\}| < 4 |\{y: {1 \over 4b_i!}C_Sx^{\alpha_{ij} - b_iN_i}y^{b_i}| < \epsilon\}| 
\eqno (5.19)$$
Integrating in $x$, we obtain
$$|\{(x,y) \in U_{ij}: |R_i(x,y)| < \epsilon\}| < 4 |\{(x,y) \in U_{ij}: {1 \over 4b_i!}C_Sx^{\alpha_{ij} 
- b_iN_i}y^{b_i}| < \epsilon\}| \eqno (5.20)$$
We bound the right-hand side of $(5.20)$ using Lemma 2.3. Since $\alpha_{ij} = a_i + N_ib_i$ and 
$b_i < a_i$, we have $b_i < \alpha_{ij} - b_iN_i$ and part c) of Lemma 2.3 applies. We obtain that
$$|\{(x,y) \in U_{ij}: |R_i(x,y)| < \epsilon\}| < C_S'\epsilon^{N_i + 1  \over \alpha_{ij}} \eqno 
(5.21)$$
Since $\alpha_{ij} = a_i + N_ib_i$, by Lemma 4.1a) the exponent ${N_i + 1  \over \alpha_{ij}}$ is at
most ${1 \over d}$, and $(5.21)$ gives us the estimate we need. Thus we are done in the setting of
Theorem 1.1.

Suppose now that we are in the setting of Theorem 1.2. Since $(a_i,b_i) \in e_{ij}$ and $S_i$ is in 
superadapted coordinates, any zero of $(S_i)_{e_{ij}}(1,y)$ is of order less than $d$. As a result, no 
matter
what $f_i$ is, there are at most finitely many $t$ for which $(R_i)_{e_{ij}}(1,y) = (S_i)_{e_{ij}}(1,y)
+ t(f_i)_{e_{ij}}(1,y)$ has a zero of order $d$ or greater on $[0,A_i]$. (This can be proven 
by an elementary argument.) Hence excluding those $t$, for a given $t$ we can divide $[0,N_i]$ into 
closed intervals $B_1,...,B_m$ such that on each $B_k$, $\partial_y^l((R_i)_{e_{ij}}(1,y))$ is nonvanishing 
for some $0 \leq l < d$. We then apply the above argument for each $B_k$, replacing $b_i$ by $l$. Thus for
each $k$ the corresponding set of points where $|R_i(x,y)| < \epsilon$ has measure less than 
$C|\epsilon|^{1 \over d}$. Adding over all $k$ we get the upper bounds required by Theorem 1.2. Thus we 
have proven the desired bounds for $M_{R_i, U_{ij}}(\epsilon)$ in the case that the bisectrix intersects 
$N(S_i)$ in the interior of a compact edge.

We now turn to the case where the bisectrix intersects $N(S_i)$ at a vertex $(d,d)$ with $(a_i,b_i) = 
(d,d)$. For fixed $x$, the $y$-cross-section of $U_{ij}$ is contained in $[c_{ij}x^{m_{ij}}, C_{ij}
x^{m_{ij}}]$ for some $c_{ij}$ 
and $C_{ij}$ which depend on the function $R_i(x,y)$. We write $[c_{ij}x^{m_{ij}}, C_{ij}x^{m_{ij}}]$
as the union of $[c_{ij}x^{m_{ij}}, x^{m_{ij}}]$ and $[x^{m_{ij}}, C_{ij}x^{m_{ij}}]$, and correspondingly
write $U_{ij} = U_{ij}^1 \cup U_{ij}^2$. We focus our attention on $U_{ij}^1$ only, as $U_{ij}^2$ is
done analogously with the roles of the two axes reversed. One technical point here is worth mentioning. 
Since $(a_i,b_i) = (d,d)$ and $S_0(x,y)$ was in superadapted coordinates, the algorithm of 
Theorem 3.1 is such that $\phi_i(x,y)$ is the identity; $F_i$ is carved out of the original disk and there
is no coordinate change. This is relevant here because it implies fractional powers of $x$ do not appear;
one can switch the $x$ and $y$ axes without any issues caused by fractional powers arising.

The argument is basically the same as that of $(5.13)-(5.21)$ so we will be brief. This time we write
$$R_i(x,y) = \sum_{a + m_{ij}b = \alpha_{ij}} r_{ab}x^ay^b + \sum_{a + m_{ij}b > \alpha_{ij},\,\,a,b < M}
r_{ab}x^ay^b + O(|x|^M + |y|^M) \eqno (5.22a)$$
$$S_i(x,y) = \sum_{a + m_{ij}b = \alpha_{ij}} s_{ab}x^ay^b + \sum_{a + m_{ij}b > \alpha_{ij},\,\,a,b < M}
s_{ab}x^ay^b + O(|x|^M + |y|^M) \eqno (5.22b)$$
The analogue to $(5.15a)-(5.15b)$ is 
$$(\partial_y^d R_i)(x,x^{m_{ij}}y) = x^{\alpha_{ij} - dm_{ij}}\partial_y^{d}((R_i)_{e_{ij}}(1,y)) + 
x^{\alpha_{ij} - dm_{ij} + \zeta}\tilde{P}(x,y) + O(|x|^{M'}) \eqno (5.23a)$$
$$(\partial_y^d S_i)(x,x^{m_{ij}}y) = x^{\alpha_{ij} - dm_{ij}}\partial_y^{d}((S_i)_{e_{ij}}(1,y)) + 
x^{\alpha_{ij} - dm_{ij} + \zeta}\tilde{Q}(x,y) + O(|x|^{M'}) \eqno (5.23b)$$
Note that since $(d,d)$ is on $e_{ij}$ with equation $a + m_{ij}b = \alpha_{ij}$, the
exponent $\alpha_{ij} - dm_{ij}$ in $(5.23a)-(5.23b)$ is equal to $d$. Hence in the setting of Theorem 
1.1, if the $\delta$ of that theorem is sufficiently small, the analogue to $(5.18b)$ here becomes
$$(\partial_y^d R_i)(x,y) >  {1 \over 4}C_S''x^d \eqno (5.24)$$
Using Lemma 2.2, the analogue to $(5.20)$ here is
$$|\{(x,y) \in U_{ij}: |R_i(x,y)| < \epsilon\}| < 4 |\{(x,y) \in U_{ij}: {1 \over 4d!}C_S''
x^dy^d| < \epsilon\}| \eqno (5.25)$$
Hence using Lemma 2.3b) now, we get that 
$$|\{(x,y) \in U_{ij}: |R_i(x,y)| < \epsilon\}| < C_S'''\epsilon^{1 \over d}|\ln \epsilon| \eqno (5.26)$$
This is the desired estimate in the setting of Theorem 1.1. As for Theorem 1.2, we make
modifications analogous to the ones before. Specifically, we again can exclude finitely many values of $t$
for which $(R_i)_{e_{ij}}(1,y)$ has zeroes of order greater than $d$ on $[c_{ij},1]$ and assume 
$(R_i)_{e_{ij}}(1,y)$ has no zeroes of order greater than $d$. One then proceeds as in the 
previous case, dividing $[c_{ij},1]$ into intervals on which either $(R_i)_{e_{ij}}(1,y)$ is nonvanishing 
or $\partial_y^k((R_i)_{e_{ij}}(1,y))$ is nonvanishing for some $k \leq d$. This completes the arguments for 
the $M_{R_i, U_{ij}}(\epsilon)$ in the exceptional cases. 

\noindent {\bf The analysis of $M_{R_i, V_{ijk}^1}(\epsilon)$.}

\noindent Analogous to $(5.14a)$, for every $i$ and $j$ we have that
$$R_i(x,x^{m_{ij}}y) = x^{\alpha_{ij}}(R_i)_{e_{ij}}(1,y) + x^{\alpha_{ij}+\zeta}P(x,y) + O(|x|^{M'})$$
There is a zero $r$ of $(R_i)_{e_{ij}}(1,y)$ such that on  $V_{ijk}^1$ we have $|y - rx^{m_{ij}}| > 
x^{m_{ij} + \eta}$. Hence in the above equation, $|y - r| > x^{\eta}$. For any $\eta' > 0$, we can choose
$\eta$ so that $|y - r| > x^{\eta}$ implies that $(R_i)_{e_{ij}}(1,y) > Cx^{\eta'}$ for the $(x,y)$ being
considered. Hence for any $\eta' > 0$ we can assume that
$$x^{\alpha_{ij}}(R_i)_{e_{ij}}(1,y) > Cx^{\alpha_{ij} + \eta'}$$
As long as $\eta' < \zeta$ and $M'$ is sufficiently large, then we therefore have
$$|R_i(x,x^{m_{ij}}y)| > C'x^{\alpha_{ij} + \eta'} \eqno (5.27)$$
Since $V_{ijk}^1$ is between the curves $y = c_1x^{m_{ij}}$ and $c_2x^{m_{ij}}$ for some $c_1$ and $c_2$,
we conclude that
$$|\{(x,y) \in V_{ijk}^1: |R_i(x,y)| < \epsilon\}| < C''\epsilon^{1 + m_{ij} \over \alpha_{ij} + \eta'}
\eqno (5.28)$$
Thus for any $(i,j)$ for which the exponent ${1 + m_{ij} \over \alpha_{ij}}$ is larger than
what we need, by shrinking $\eta'$ enough, $(5.28)$ gives that $M_{R_i,V_{ijk}^1}(\epsilon)$ satisfies a better
estimate than what we need. However, we saw that the only situations in which this exponent is not
better than what we need are the exceptional cases discussed above. But these are exactly the cases in which
we proved the desired upper bounds for all of $M_{R_i, U_{ij}}(\epsilon)$. Hence by shrinking $\eta$ enough,
$(5.28)$ gives better than the needed estimates for all remaining $V_{ijk}^1$ and we are done.

\noindent {\bf The analysis of $M_{R_i, V_{ijk}^2}(\epsilon)$.}

\noindent We now focus our attention on some fixed $V_{ijk}^2$, which consists of the points in 
$U_{ij}$ for which $|y - rx^{m_{ij}}| < x^{m_{ij} + \eta}$ for some $r \in \R$. The exact value of $r$ 
will not be important in what follows. We apply the resolution of 
singularities algorithm of Theorem 3.1 to $f_i(x,y)$ on $F_i$, 
and consider the sets called $D_i$ in that theorem that intersect $V_{ijk}$. Since we are already 
using the index $i$, we refer to them as $D_l$ here. To each $D_l$ there is a coordinate change $\phi_l$ 
such that $f_i \circ  \phi_l$ is comparable to a monomial on $D_l$ in the sense of Theorem 3.1 b) or c). 
Since we consider only those $D_l$ intersecting $V_{ijk}$, the function $\phi_l$ is such that $|y| < 
x^{m_{ij} + \eta}$ on $\phi_l^{-1}D_l$. Write $S_i'(x,y) = S_i \circ \phi_l(x,y)$. Our first task will be to 
understand $S_i'(x,y)$'s behavior on the set $\phi_l^{-1}D_l$, which we denote by $E_l$. 

To this end, note that the domain $F_i$ of $S_i(x,y)$ has upper boundary given by $y = A_i x^{N_i} +$
higher order terms, and lower boundary given by either $y = 0$ or $y = a_ix^{n_i} +$ higher order 
terms, where $n_i > N_i$. Since $S_i(x,y) \sim x^{a_i}y^{b_i}$ on $F_i$, if $n_i > m_{ij} > N_i$ then on
$N(S_i)$ the linear function $a + m_{ij}b$ is minimized at exactly one point, $(a_i,b_i)$. If $m_{ij} = N_i$
or $n_i$, then $a + m_{ij}b$ will be minimized at $(a_i,b_i)$ and possibly also at other points in $N(S_i)$.

Next, suppose $f(x)$ is any function such that $f(x^K)$ is smooth for some $K$ and such that the Taylor 
expansion of $f(x)$ has initial term $cx^{m_{ij}}$ for some $c$. 
Also suppose that $n_i > m_{ij} > N_i$. Then the Newton polygon of $S_i(x,y - f(x))$ will have 
an edge of slope $-{1 \over m_{ij}}$ containing the point $(a_i + m_{ij}b_i, 0)$ and no edges with 
lesser slope (i.e. no more horizontal edges). This fact implies that $S_i(x,y - f(x)) \sim
x^{a_i + m_{ij}b_i}$ on $V_{ijk}^2$, since $|y| < |x|^{m_{ij} + \eta}$ on $V_{ijk}^2$. Note that the same 
is true if $m_{ij} = N_i$, if $c$ is small enough, or if $m_{ij} = n_i$,
if $c$ is large enough. For in these cases the vertices other than $(a_i,b_i)$ minimizing $a + m_{ij}b$
do not interfere. Consequently, we may assume that  $S_i'(x,y) = S_i \circ \phi_l(x,y)$ and its Newton 
polygon has the above properties. For $S_i' = S_i \circ \phi_l(x,y)$ is either of the form 
$S_i(x,y - f(x))$ or $S_i(x,-y - f(x))$ for an $f(x)$ of this type. (In the case where $m_{ij} = n_i$ or
$N_i$, if the quantity called $\xi$ in the proof
of the resolution of singularities algorithm that was originally applied to $S_i(x,y)$ was small 
enough, which we may assume, then the coefficient $a_i$ will be large enough and the coefficient 
$A_i$ will be small enough for this to work.) 

Let $f_i' = f_i \circ \phi_l$ and $R_i' = R_i \circ \phi_l$, so that $R_i' = S_i' + t f_i'$. We now will 
estimate the various $M_{R_i', E_l}(\epsilon)$ and show they satisfy estimates better than those that we 
need. We start with $E_l$ satisfying part c) of Theorem 3.1; that is, when $f_i'(x,y) \sim x^{\alpha_l}
y^{\beta_l}$ with $\beta_l > 0$. Define $X_l$ to be the set of points in $E_l$ for which $|R_i'(x,y)| > 
{1 \over 2}|S_i'(x,y)|$, and let $Y_l = E_l - X_l$. Then we have
$$|\{(x,y) \in X_l: |R_i'(x,y)| < \epsilon\}| < |\{(x,y) \in X_l: |S_i'(x,y)| < 2\epsilon\}| \eqno (5.29)$$
$$\leq M_{S_i',E_l}(\epsilon) $$
This is better than the estimate we need, so we focus our attention on $Y_l$. Note that on $Y_l$ we have
$${1 \over 2}|S_i'(x,y)| \leq |tf_i'(x,y)| \leq {3 \over 2}|S_i'(x,y)| \eqno (5.30)$$
Part c) of Theorem 3.1 says that, modulo a function vanishing to infinite order at $(0,0)$, $f_i'(x,y) \sim x^{\alpha_l}y^{\beta_l}$ on $E_l$, while by the above 
discussion $S_i'(x,y) \sim x^{a_i + m_{ij}b_i}$ on $E_l$. Thus when $(5.30)$ holds we have 
$x^{\alpha_l}y^{\beta_l} \sim x^{a_i + m_{ij}b_i}$. Next, note that by Theorem 3.1c), for any $K$ one has
$$|\partial_y f_i'(x,y)| > C_0x^{\alpha_l}y^{\beta_l-1} - O(x^K) > C_1x^{a_i + m_{ij}b_i}y^{-1} - O(x^K)$$
$$> C_2x^{a_i + m_{ij}b_i - m_{ij} - \eta} \eqno (5.31)$$
The last inequality follows from the fact that $|y| < |x|^{m_{ij} + \eta}$ on $E_l$.
On the other hand, since $S_i'(x,y)$'s Newton polygon has an edge of slope $-{1 \over m_{ij}}$ containing
$(a_i + m_{ij}b_i,0)$, and no other edges more horizontal that this one, we also have
$$|\partial_y S_i'(x,y)| < C_3x^{a_i + m_{ij}b_i - m_{ij}} \eqno (5.32)$$
Thus in taking the $y$ derivative of $R_i'(x,y) = S_i'(x,y) + tf_i'(x,y)$, the derivative of the second term
dominates and (for $t \neq 0$) we have
$$|\partial_y R_i'(x,y)| > C_4 x^{a_i + m_{ij}b_i - m_{ij} - \eta} \eqno (5.33)$$
Next, we apply Theorem 2.2 and for some $x_0 > 0$ get that
$$|\{(x,y) \in Y_l: |R_i'(x,y)| < \epsilon\}| \leq 4|\{(x,y) \in Y_l: C_4x^{a_i + m_{ij}b_i - m_{ij} - \eta}y  < 
\epsilon\}|$$
$$\leq |(x,y): 0 < x < x_0,\,\,|y| < x^{m_{ij} + \eta},\,\,C_5x^{a_i + m_{ij}b_i - m_{ij} - \eta}y  < 
\epsilon\}|\eqno (5.34)$$
The last inequality again uses that $|y| < |x|^{m_{ij} + \eta}$ on $E_l$.
Equation $(5.34)$ can be bounded with the help of Lemma 2.3. If part a) or b) applies, it is bounded by
$C_6 \epsilon |\ln \epsilon|$, better than the estimate we need since we are assuming $d > 1$. If part c) 
applies, we get 
$$|\{(x,y) \in Y_l: |R_i'(x,y)| < \epsilon\}| < C_7 \epsilon^{m_{ij} + \eta + 1 \over a_i + m_{ij}b_i}
\eqno (5.35)$$
Since $\eta > 0$, this exponent is better than ${m_{ij} + 1 \over a_i + m_{ij}b_i}$, the reciprocal of the
ordinate of the intersection of the bisectrix with the line of slope $-{1 \over m_{ij}}$ containing 
$(a_i,b_i)$. By $(5.11)$, this ordinate is at most $d$. Hence the exponent on the right-hand side of 
$(5.35)$ is greater than ${1 \over d}$, better than the estimate that we need.

We now move to the case where part b) of Theorem 3.1 is satisfied; that is, we assume $\beta_l = 0$ and thus
$f_i'(x,y) \sim x^{\alpha_l}$ on $E_l$. If $\alpha_l$ were less than $a_i + m_{ij}b_i$, for $t \neq 0$ we would have
$R_i'(x,y) = S_i'(x,y) + t f_i'(x,y) \sim x^{\alpha_l}$ on $E_l$ as well. Thus for small enough $x$ (which we
may assume by making the radius of the original disk $D$ sufficiently small) we would have $|R_i'(x,y)| 
\geq |S_i'(x,y)|$, and thus $M_{R_i',E_l}(\epsilon) \leq M_{S_i',E_l}(\epsilon)$, better than the estimate that we need. 
Similarly, if $\alpha_l$ were greater than $a_i + m_{ij}b_i$, $f_i'(x)$ would be small compared to $S_i'(x,y)$
and thus when $x$ is sufficiently small we have $|R_i'(x,y)| \geq {1 \over 2}|S_i'(x,y)|$, once again 
giving an estimate better than the one we need. Hence in the following we assume that $\alpha_l = a_i +
m_{ij}b_i$.

Next, note that since we are in the setting of part b) of Theorem 3.1, the upper boundary of $E_l$ has 
equation $y = b x^p +$ higher order terms, where $p \geq m_{ij} + \eta$. The lower boundary of $E_l$ is the
$x$-axis. Since $f_i'(x,y) \sim x^{a_i + m_{ij}b_i}$ on all of $E_l$, the Newton polygon of $f_i'$ contains
the vertex $(a_i + m_{ij}b_i,0)$, and has no edge of slope greater than $-{1 \over p}$. 

We now look at $S_i'(x,x^py)$ and $f_i'(x,x^py)$. First, in analogy with $(5.13a)-(5.13b)$ we write
$$f_i'(x,y) = \sum_{a + pb = a_i + m_{ij}b_i} f_{ab}'x^ay^b + \sum_{a + pb > a_i + m_{ij}b_i,\,\,a,b < M}
f_{ab}'x^ay^b + O(|x|^M + |y|^M) \eqno (5.36a)$$
$$S_i'(x,y) = \sum_{a + pb = a_i + m_{ij}b_i} s_{ab}'x^ay^b + \sum_{a + pb > a_i + m_{ij}b_i,\,\,a,b < M}
s_{ab}'x^ay^b + O(|x|^M + |y|^M) \eqno (5.36b)$$
It is worth pointing out that since $N(S_i')$ has no edge of slope greater than $-{1 \over m_{ij}}$, the 
first sum in $(5.36b)$ has only one term. Next, in analogy with $(5.14a)-(5.14b)$, we have
$$f_i'(x,x^py) = x^{a_i + m_{ij}b_i}(f_i')_{e_p}(1,y) + x^{a_i + m_{ij}b_i+\zeta}P'(x,y) + O(|x|^{M'})
\eqno (5.37a)$$
$$S_i'(x,x^py) = x^{a_i + m_{ij}b_i}(S_i')_{e_p}(1,y) + x^{a_i + m_{ij}b_i+\zeta}Q'(x,y) + O(|x|^{M'})
\eqno (5.37b)$$
Here $(f_i')_{e_p}(x,y)$ denotes the sum of all terms $f_{ab}'x^ay^b$ of the Taylor expansion of $f_i'(x,y)$
at the origin with $a + pb = a_i + m_{ij}b_i$ (the minimal possible value of $a + pb$), and similarly for
$(S_i')_{e_p}(x,y)$. Note that $(S_i')_{e_p}(1,y)$ is constant here. Hence for all but finitely many $t$,
$(S_i')_{e_p}(1,y) + t(f_i')_{e_p}(1,y)$ has no zeroes of order greater than one (this is actually true in
more general scenarios as well.) Hence excluding these values of $t$, we write
$$R_i'(x,x^py) = x^{a_i + m_{ij}b_i}(R_i')_{e_p}(1,y) + x^{a_i + m_{ij}b_i+\zeta}R'(x,y) + O(|x|^{M'})
\eqno (5.38)$$
Here $(R_i')_{e_p}(1,y)$ has no zeroes of order greater than one. The $y$-range in $(5.38)$ is contained
in $[0,b + \epsilon]$ for any $\epsilon > 0$. (Recall $b$ is such that the upper boundary of $E_l$ is given 
by $y = b x^p 
+...$). We write $[0,b + \epsilon]$ as a union $B_1 \cup ... \cup B_k$, where on each $B_k$ either 
$(R_i')_{e_p}(1,y)$ is nonvanishing, or $\partial_y(R_i')_{e_p}(1,y)$ is nonvanishing. 

Consider the case of a $B_k$ for which $(R_i')_{e_p}(1,y)$ is nonvanishing. Then on the domain of $(5.38)$,
we have $R_i'(x,x^py) > C_1x^{a_i + m_{ij}b_i}$. Translating back into the coordinates of $E_l$, we have
$R_i'(x,y) > C_1x^{a_i + m_{ij}b_i}$ on a subset $A_k$ of $\{(x,y): 0 < x < x_0, 0 < y < (b + 1)x^p\}$ for
some $x_0 > 0$. Thus we have
$$M_{R_i',A_k}(\epsilon) = |\{(x,y) \in A_k: |R_i'(x,y)| < \epsilon\}| < C_2\epsilon^{1 + p \over a_i + m_{ij}b_i} 
\eqno (5.39)$$
Since $p > m_{ij}$, this exponent is better than ${1 + m_{ij} \over a_i + m_{ij}b_i} \geq {1 \over d}$.
Thus $(5.39)$ is better than the exponent we need. Now consider the case of a $B_k$ for which $\partial_y
((R_i')_{e_p}(1,y))$ is nonvanishing. In this case, the relevant analogue of $(5.38)$ is
$$(\partial_y R_i')(x,x^py) = x^{a_i + m_{ij}b_i - p}\partial_y ((R_i')_{e_p}(1,y)) + x^{a_i + m_{ij}b_i - p +\zeta}R''(x,y) 
+ O(|x|^{M'}) \eqno (5.40)$$
Hence on the domain of $(5.40)$ we have $|(\partial_y R_i')(x,x^py)| > C_3x^{a_i + m_{ij}b_i - p}$ if 
$x$ is sufficiently small, which we may assume. Translating
back into the original coordinates of $E_l$, this time we get that $\partial_y R_i'(x,y) > 
C_1x^{a_i + m_{ij}b_i - p}$ on a subset $A_k$ of $\{0 < x < x_0, 0 < y < (b + 1)x^p\}$. Thus we may apply
Lemma 2.2 and say that
$$M_{R_i',A_k}(\epsilon) = |\{(x,y) \in A_k: |R_i'(x,y)| < \epsilon\}| < 4|\{(x,y) \in A_k: C_4x^{a_i + 
m_{ij}b_i - p}y < \epsilon\}|$$
$$\leq 4|\{(x,y) \in E_l: C_4x^{a_i + m_{ij}b_i - p}y < \epsilon\}| \eqno (5.41)$$
Since $E_l$ is a subset of $\{(x,y): 0 < x < x_0,\,\,0 < y < (b + 1)x^p\}$ for some $x_0 > 0$, we may use 
Lemma 2.3 to estimate the right-hand side of $(5.41)$. If part a) or b) of the lemma applies, we get that 
$|M_{R_i',A_k}(\epsilon)| < C\epsilon|\ln \epsilon|$, better than the estimate we need. If part c) applies, 
we get 
$$M_{R_i',A_k}(\epsilon) < C_5 \epsilon^{1 + p \over a_i + m_{ij}b_i} \eqno (5.42)$$
This exponent is the same as that of $(5.39)$, which we saw is better than what we need. This completes
the proofs of Theorem 1.1 and 1.2. 

\noindent {\bf 6. References.}

\noindent [AGV] V. Arnold, S Gusein-Zade, A Varchenko, {\it Singularities of differentiable maps},
Volume II, Birkhauser, Basel, 1988. \parskip = 3pt\baselineskip = 3pt

\noindent [CaCWr] A. Carbery, M. Christ, J. Wright, {\it Multidimensional van der Corput and sublevel set 
estimates.}, J. Amer. Math. Soc. {\bf 12} (1999), no. 4, 981-1015.

\noindent [G1] M. Greenblatt, {\it The asymptotic behavior of degenerate oscillatory integrals in two
dimensions}, to appear, J. Funct. Anal. 

\noindent [G2] M. Greenblatt, {\it A direct resolution of singularities for functions of two variables 
with applications to analysis}, J. Anal. Math. {\bf 92} (2004), 233-257.

\noindent [G3] M. Greenblatt, {\it Sharp $L\sp 2$ estimates for one-dimensional oscillatory integral 
operators with $C\sp \infty$ phase.} Amer. J. Math. {\bf 127} (2005), no. 3, 659-695. 

\noindent [H] L. H$\ddot{{\rm o}}$rmander, {\it The analysis of linear partial differential operators. I. 
Distribution theory and Fourier analysis}, 2nd ed. Springer-Verlag, Berlin, (1990).
xii+440 pp. 

\noindent [IKeM] I. Ikromov, M. Kempe, and D. M\"uller, {\it Sharp $L^p$ estimates for maximal operators
associated to hypersurfaces in $\R^3$ for $p > 2$}, preprint.

\noindent [IoSa] A. Iosevich, E. Sawyer, {\it Oscillatory integrals and maximal averages over homogeneous
surfaces}, Duke Math. J. {\bf 82} no. 1 (1996), 103-141. 

\noindent [K1] V. N. Karpushkin, {\it A theorem concerning uniform estimates of oscillatory integrals when
the phase is a function of two variables}, J. Soviet Math. {\bf 35} (1986), 2809-2826.

\noindent [K2] V. N. Karpushkin, {\it Uniform estimates of oscillatory integrals with parabolic or 
hyperbolic phases}, J. Soviet Math. {\bf 33} (1986), 1159-1188.

\noindent [K3] V. N. Karpushkin, {\it Uniform estimates for volumes}, Tr. Math. Inst. Steklova {\bf 221} 
(1998), 225-231.

\noindent [PS] D. H. Phong, E. M. Stein, {\it The Newton polyhedron and oscillatory integral operators},
Acta Mathematica {\bf 179} (1997), 107-152.

\noindent [PSSt] D. H. Phong, E. M. Stein, J. Sturm, {\it On the growth and stability of real-analytic 
functions}, Amer. J. Math. {\bf 121} (1999), no. 3, 519-554.

\noindent [R] K. M. Rogers, {\it Sharp Van der Corput estimates and minimal divided differences}, 
Proc. Amer. Math. Soc. {\bf 133} (2005), no. 12, 3543-3550.

\noindent [V] A. N. Varchenko, {\it Newton polyhedra and estimates of oscillatory integrals}, Functional 
Anal. Appl. {\bf 18} (1976), no. 3, 175-196.

\line{}
\line{}

\noindent Department of Mathematics, Statistics, and Computer Science \hfill \break
\noindent University of Illinois at Chicago \hfill \break
\noindent 322 Science and Engineering Offices \hfill \break
\noindent 851 S. Morgan Street \hfill \break
\noindent Chicago, IL 60607-7045 \hfill \break
\end